\documentclass[12pt]{article}
\usepackage{amsmath,amssymb,amscd}

\newcommand{\eqdef}{\stackrel{\text{def}}{=}}
\newcommand{\n}{\nonumber \\}
\newcommand{\bm}{\boldsymbol}
\newcommand{\ignore}[1]{}
\renewcommand{\theequation}{\arabic{section}.\arabic{equation}}

\allowdisplaybreaks[4]

\begin{document}

\baselineskip=20pt

\newfont{\elevenmib}{cmmib10 scaled\magstep1}
\newcommand{\preprint}{
     \begin{flushleft}
       \elevenmib Yukawa\, Institute\, Kyoto\\
     \end{flushleft}\vspace{-1.3cm}
     \begin{flushright}\normalsize  \sf
       DPSU-07-05\\
       YITP-07-91\\
       December 2007
     \end{flushright}}
\newcommand{\Title}[1]{{\baselineskip=26pt
     \begin{center} \Large \bf #1 \\ \ \\ \end{center}}}
\newcommand{\Author}{\begin{center}
     \large \bf Satoru Odake${}^a$ and Ryu Sasaki${}^b$ \end{center}}
\newcommand{\Address}{\begin{center}
       $^a$ Department of Physics, Shinshu University,\\
       Matsumoto 390-8621, Japan\\
       ${}^b$ Yukawa Institute for Theoretical Physics,\\
       Kyoto University, Kyoto 606-8502, Japan
     \end{center}}
\newcommand{\Accepted}[1]{\begin{center}
     {\large \sf #1}\\ \vspace{1mm}{\small \sf Accepted for Publication}
     \end{center}}

\preprint
\thispagestyle{empty}
\bigskip\bigskip\bigskip

\Title{Orthogonal Polynomials from Hermitian Matrices}
\Author

\Address
\vspace{1cm}

\begin{abstract}
A unified theory of orthogonal polynomials of a discrete variable
is presented through the eigenvalue problem of hermitian matrices
of finite or infinite dimensions.
It can be considered as a matrix version of exactly solvable
Schr\"odinger equations. The hermitian matrices (factorisable
Hamiltonians) are real symmetric tri-diagonal (Jacobi) matrices
corresponding to second order difference equations.
By solving the eigenvalue problem in two different ways, the duality
relation of the eigenpolynomials and their dual polynomials is
explicitly established.
Through the techniques of exact Heisenberg operator solution and
shape invariance, various quantities, the two types of eigenvalues
(the eigenvalues and the sinusoidal coordinates), the coefficients
of the three term recurrence, the normalisation measures and the
normalisation constants etc. are determined explicitly.
\end{abstract}

\section{Introduction}
\label{intro}

Due to the long history and rich and diverse applications, there
are many different ways to introduce orthogonal polynomials \cite{szego}.
Among many types and kinds of orthogonal polynomials, we will focus
in this paper on the so-called orthogonal polynomials of a discrete
variable \cite{nikiforov,askey,gasper,ismail}, whose orthogonality
measures are concentrated on discrete points, either finite or
infinite in number.
We will present a unified theory of these polynomials based on the
eigenvalue problem \eqref{eigenpr} of a special class of hermitian
matrices, which are real symmetric {\em tri-diagonal\/} (Jacobi)
matrices \eqref{genham}, \eqref{Jacobiform}.
Since the spectrum of a Jacobi matrix is simple, the orthogonality
of eigenvectors is guaranteed. The eigenvalue problem of a Jacobi
matrix can be considered as a difference equation version of the
Schr\"odinger equation for one degree of freedom system, the most
basic equation of quantum mechanics \cite{flugge}.
Various classical orthogonal polynomials (the Hermite, Laguerre,
Jacobi and their restrictions) have appeared as the eigenfunctions
of exactly solvable quantum mechanics \cite{flugge}.

Roughly speaking, our line of arguments is a deformation (discretisation)
of the main trend of the twentieth century mathematical physics;
to pursue the parallelism between the matrix eigenvalue problem and
the ordinary differential equations of the Sturm-Liouville type,
with {\em differential\/} equations replaced by {\em difference\/}
equations.

We apply many ideas and methods for solving Schr\"odinger equations:
Crum's theorem \cite{crum}, the factorisation method \cite{infhul}
(or the so-called supersymmetric quantum mechanics \cite{susyqm}),
the method of exact Heisenberg operator solutions and the
creation/annihilation operators \cite{os7,os9} together with
symmetries (shape invariance \eqref{shapeinv1} \cite{genden,os456}
and closure relation \eqref{closurerel1} \cite{os7}) to elucidate
the universal structure of the orthogonal polynomials of a discrete
variable.
The Jacobi matrix eigenvalue problem
$\mathcal{H}\psi=\mathcal{E}\psi$ \eqref{eigenpr} can be solved in
a different way, through the three term recurrence relations for
the dual polynomial in $\mathcal{E}$.
Combining the {\em duality\/}, which can be stated at many different
levels, with the above mentioned solution techniques, various
quantities of the orthogonal polynomials, the coefficients of the
three term recurrence, the orthogonality measure and the normalisation
constants, the difference equations for the dual polynomials, etc.
are given explicitly in an elementary manner.
We stress that duality is established equally for the finite and
infinite dimensional cases.
All the examples discussed in this paper are taken from the review of
Koekoek and Swarttouw \cite{koeswart} and they are known to satisfy
the closure relation \eqref{closurerel1}.
The orthogonal polynomials with the Jackson integral type measure,
{\em e.g.} the Big $q$-Jacobi polynomial, etc. will need different
treatment and will be discussed elsewhere.
As an inverse step to characterise the orthogonal polynomials of a
discrete variable, we solve the closure relation algebraically to
determine the possible forms of the Jacobi matrices.

It should be stressed that the present duality has far richer contents
than the well-known duality of the column and row eigenvectors.
That is, in the (hermitian) eigenvalues problems, the orthogonality of
the complete set of column eigenvectors automatically implies that of
the row eigenvectors.
In the present case, the column eigenvectors are polynomials in
{\em sinusoidal coordinate\/} \eqref{lowtri} $\eta(x)$ and the row
eigenvectors are polynomials in the {\em eigenvalues\/} \eqref{lowtri}
$\mathcal{E}(n)$, and the duality is formulated as the equality of
these two polynomials on the integer lattice points \eqref{Duality}.

The formulation of the unified theory of orthogonal polynomials based on
the eigenvalue problems of Jacobi matrices is new.
All the quantities and formulas are derived from the two input functions
$B(x)$ and $D(x)$ \eqref{genham}, \eqref{BDcondition}.
However, a good part of the explicit formulas for specific polynomials
explored in section five is already known individually in the existing
theories of orthogonal polynomials based on various settings.

This paper is organised as follows. In section two, the general form
of the hermitian matrices (or the Hamiltonians) is given and the
solution procedure to arrive at the orthogonal polynomials is
outlined. They are polynomials in the {\em sinusoidal coordinate\/}
\eqref{lowtri} which plays an important role in exactly solvable
quantum mechanics \cite{os7,os9}.
The orthogonality measure is given by the solution $\phi_0(x)$ of
the factorised equation \eqref{Aphi0=0}. The $\phi_0$ is the ground
state wavefunction.
In section three the {\em dual polynomials\/} (in the eigenvalue
$\mathcal{E}(n)$) are introduced by a different solution method
(the three term recurrence \eqref{dual3term}) of the same eigenvalue
problem \eqref{eigenpr}. The duality at many different levels,
{\em e.g.} $x\leftrightarrow n$, $\eta(x)\leftrightarrow \mathcal{E}(n)$,
is displayed \eqref{dualityxn}--\eqref{dualityphi0dn}.
Section four explores the underlying symmetry properties of the
Hamiltonian, shape invariance \S\ref{shapesection}, closure relation
\S\ref{closuresection}. The general formulas of the coefficients $A_n$
and $C_n$ of the three term recurrence are derived.
The forms of the sinusoidal coordinate and the diagonal elements of
the Hamiltonian are determined in \S\ref{detetaBD} from the closure
relation.
The dual closure relation is discussed in \S\ref{dualclosuresection}.
Section five provides the explicit forms of various quantities
derived in previous sections for the concrete examples of the
orthogonal polynomials in the orders given in the review \cite{koeswart}.
Some orthogonal polynomials not listed in the review \cite{koeswart}
are also discussed in some detail in \S\ref{alternatives}.
Appendix A gives the possible forms of the Jacobi matrices
(the Hamiltonians) \eqref{Bxformula2}--\eqref{Dxformula2}
as the solutions of the closure relation.
The results show the boundedness and unboundedness of the Hamiltonians
and eigenvalues clearly.
Appendix B provides the collection of the definitions of basic symbols
and functions for self-containedness.
Throughout this paper we use the parameter $q$ in the range $0<q<1$.

\section{Hamiltonian}
\label{Hamil}
\setcounter{equation}{0}

The starting point is the eigenvalue problem for a hermitian matrix
$\mathcal{H}$
\begin{equation}
  \mathcal{H}\psi=\mathcal{E}\psi.
  \label{eigenpr}
\end{equation}
The rows and columns of $\mathcal{H}$ are indexed by non-negative
integers $x$ and $y$, either finite
\begin{equation}
  \text{case (I)}\quad x, y=0,1,2,\ldots,N,
\end{equation}
or infinite:
\begin{equation}
  \text{case (II)}\quad x, y=0,1,2,\ldots.
\end{equation}
Let us call the hermitian matrix $\mathcal{H}$ {\em Hamiltonian\/},
since the eigenvalue problem \eqref{eigenpr} can be considered as
the Schr\"odinger equation which is a {\em difference} equation,
instead of differential in ordinary quantum mechanics.
The Hamiltonian we consider in this paper has a general form
\begin{equation}
  \mathcal{H}\eqdef
  -\sqrt{B(x)}\,e^{\partial}\sqrt{D(x)}
  -\sqrt{D(x)}\,e^{-\partial}\sqrt{B(x)}
  +B(x)+D(x),
  \label{genham}
\end{equation}
in which the two functions $B(x)$ and $D(x)$ are real and {\em positive}
but vanish at the boundary:
\begin{align}
  B(x)>0,\quad D(x)>0,\quad  D(0)=0\ ;\quad
  B(N)=0\ \ \text{for case (I)}.
  \label{BDcondition}
\end{align}
In \eqref{genham} $\partial\eqdef \frac{d}{dx}$ is the differential
operator and its exponentiation gives finite shifts:
\[
  e^{\pm\partial}\psi(x)=\psi(x\pm1).
\]
Thus, in fact, the eigenvalue problem \eqref{eigenpr} can be written as
a difference equation on non-negative integer points:
\begin{align}
  \bigl(B(x)+D(x)\bigr)\psi(x)-\sqrt{B(x)D(x+1)}\,\psi(x+1)
  &-\sqrt{B(x-1)D(x)}\,\psi(x-1)
  =\mathcal{E}\psi(x),\n
  &x=0,1,\ldots,(N),\ldots.
  \label{diffeq}
\end{align}
The boundary condition $D(0)=0$ is necessary for the term $\psi(-1)$
does not appear, and $B(N)=0$ is necessary  for the term $\psi(N+1)$
does not appear in the finite dimensional matrix case.
It is also easy to see that $\psi(0)$ does not vanish for any eigenvector
\begin{equation}
  \psi(0)\neq0,\qquad \bigl(\psi(N)\neq0\ \ \text{for case (I)}\bigr).
  \label{psizero}
\end{equation}

Let us remark that  the ($q$-)Askey-scheme of hypergeometric
orthogonal polynomials with absolutely continuous weight functions,
for example, the Meixner-Pollaczek, the Wilson and the Askey-Wilson
polynomials, have already been well understood as the eigenvalue
problem of hermitian (self-adjoint) operators $\mathcal{H}$
\cite{os456,os7}.
The generic forms of the Hamiltonian $\mathcal{H}$ are:
\begin{align}
  \mathcal{H}&\eqdef\sqrt{V(x)}\,e^{-i\partial}\sqrt{V(x)^*}
  +\!\sqrt{V(x)^*}\,e^{i\partial}\sqrt{V(x)}-V(x)-V(x)^*,
  \label{H1}\\
  \mathcal{H}&\eqdef\sqrt{V(x)}\,q^{-i\partial}\!\sqrt{V(x)^*}
  +\sqrt{V(x)^*}\,q^{i\partial}\!\sqrt{V(x)}-V(x)-V(x)^*.
  \label{Hq1}
\end{align}
These also give difference Schr\"odinger equations, but the shift is
in the {\em imaginary\/} direction, instead of real in the case
\eqref{diffeq}.

Although the Hamiltonian $\mathcal{H}$ \eqref{genham} is presented in
a difference operator form, it is in fact a real symmetric
{\em tridiagonal\/} (Jacobi) matrix:
\begin{align}
  \mathcal{H}&=(\mathcal{H}_{x,y}),\qquad
  \mathcal{H}_{x,y}=\mathcal{H}_{y,x},\\
  \mathcal{H}_{x,y}&=
  -\sqrt{B(x)D(x+1)}\,\delta_{x+1,y}-\sqrt{B(x-1)D(x)}\,\delta_{x-1,y}
  +\bigl(B(x)+D(x)\bigr)\delta_{x,y}.
  \label{Jacobiform}
\end{align}
It is well-known that the spectrum of a Jacobi matrix is simple;
that is, there is no degeneracy of the eigenvalues
\[
  \mathcal{E}(0)<\mathcal{E}(1)<\mathcal{E}(2)<\cdots.
\]

\bigskip
The first step for solving the eigenvalue equation \eqref{eigenpr}
is to rewrite the Hamiltonian \eqref{genham} in a factorised form:
\begin{align}
  \mathcal{H}&=\mathcal{A}^{\dagger}\mathcal{A},
  \label{factor}\\
  \mathcal{A}^{\dagger}&=\sqrt{B(x)}-\sqrt{D(x)}\,e^{-\partial},
  \label{Ad}\\
  \mathcal{A}&=\sqrt{B(x)}-e^{\partial}\sqrt{D(x)},
  \label{A}
\end{align}
with the forward (backward) shift operator $\mathcal{A}$
($\mathcal{A}^{\dagger}$) being hermitian conjugate of each other.
Or in the matrix form, the tridiagonal (Jacobi) matrix
\eqref{Jacobiform} is decomposed into a product of
$\mathcal{A}^{\dagger}$ and $\mathcal{A}$
\begin{align}
  (\mathcal{A}^{\dagger})_{x,y}&=
  \sqrt{B(x)}\,\delta_{x,y}-\sqrt{D(x)}\,\delta_{x-1,y},\\
  \mathcal{A}_{x,y}&=
  \sqrt{B(x)}\,\delta_{x,y}-\sqrt{D(x+1)}\,\delta_{x+1,y},
\end{align}
with $\mathcal{A}^{\dagger}$ ($\mathcal{A}$) having non-vanishing
entries on the diagonal and one below (above) the diagonal.
Throughout this paper we adopt the (standard) convention that
the $(0,0)$ element of the matrix is at the upper left corner.
The factorisation also means that the Hamiltonian \eqref{genham}
is semi-positive definite. In fact, the ground state (the lowest
eigenvector) $\phi_0(x)$ is annihilated by $\mathcal{A}$ and thus
it is a zero-mode of the Hamiltonian:
\begin{equation}
  \mathcal{A}\phi_0(x)=0 \Longrightarrow \mathcal{H}\phi_0(x)=0,\qquad
  \mathcal{E}(0)=0.
  \label{Aphi0=0}
\end{equation}
As we will see shortly, $\phi_0(x)^2$ is the orthogonality measure
for the eigenpolynomials \eqref{phi0^2PnPn=dn^2}.
This is the same situation as in the ordinary quantum mechanics and
also in the above mentioned theory \eqref{H1}-\eqref{Hq1} for
the Wilson and Askey-Wilson polynomials, etc. \cite{os456,os7}.
The above equation \eqref{Aphi0=0}, being a two term recurrence relation,
\begin{equation}
  \frac{\phi_0(x+1)}{\phi_0(x)}=\sqrt{\frac{B(x)}{D(x+1)}}
  \label{phi0/phi0=B/D},
\end{equation}
can be solved elementarily with the boundary (initial) condition
\begin{align}
  \phi_0(0)&=1,
  \label{iniconphi0}\\
  \phi_0(x)&=\sqrt{\prod_{y=0}^{x-1}\frac{B(y)}{D(y+1)}},\quad
  x=1,2,\ldots.
  \label{phi0=prodB/D}
\end{align}
With the standard convention $\prod_{k=n}^{n-1}*=1$, the expression
\eqref{phi0=prodB/D} is valid for $x=0$ also.
For the infinite matrix case (II), the requirement of the finite
$\ell^2$ norm of the eigenvectors
\begin{equation}
  \sum_{x=0}^{\infty}\phi_0(x)^2=\sum_{x=0}^{\infty}
  \,\prod_{y=0}^{x-1}\frac{B(y)}{D(y+1)}<\infty
  \label{phizero2}
\end{equation}
imposes constraints on the asymptotic behaviours of $B(x)$ and $D(x)$.

Next let us determine the other eigenvectors (the excited states)
in a factored form
\begin{equation}
  \psi(x)=\phi_0(x)v(x).
\end{equation}
The eigenvalue problem \eqref{eigenpr} for $\psi$ is rewritten to that
for $v$
\begin{equation}
  \widetilde{\mathcal H}\,v(x)=\mathcal{E}\,v(x),
  \label{vequation}
\end{equation}
in which $\widetilde{\mathcal H}$ is the similarity transformed
Hamiltonian in terms of the ground state wavefunction $\phi_0(x)$:
\begin{align}
  \widetilde{\mathcal{H}}&\eqdef
  \phi_0(x)^{-1}\circ \mathcal{H}\circ\phi_0(x)\n
  &=B(x)(1-e^{\partial})+D(x)(1-e^{-\partial}).
\end{align}
It gives rise to a simple difference equation for $v(x)$
\eqref{vequation}.
Obviously a constant is a solution
\begin{equation}
  \widetilde{\mathcal{H}}\,1=0.
\end{equation}
In the matrix form $\widetilde{\mathcal{H}}$ is another tridiagonal
matrix
\begin{equation}
  \widetilde{\mathcal{H}}=(\widetilde{\mathcal{H}}_{x,y}),\quad
  \widetilde{\mathcal{H}}_{x,y}=B(x)(\delta_{x,y}-\delta_{x+1,y})+
  D(x)(\delta_{x,y}-\delta_{x-1,y}).
\end{equation}

For all the cases discussed in this paper, the similarity transformed
Hamiltonian $\widetilde{\mathcal{H}}$ is {\em lower triangular\/}
with respect to the special basis
\begin{equation}
  1,\ \eta(x),\ \eta(x)^2,\ \ldots, \eta(x)^n,\ \ldots,
\end{equation}
spanned by the {\em sinusoidal coordinate\/} $\eta(x)$
\cite{os456,os7}:
\begin{equation}
  \widetilde{\mathcal{H}}\eta(x)^n=\mathcal{E}(n)\eta(x)^n+
  \text{lower orders in}\ \eta(x).
  \label{lowtri}
\end{equation}
Here $\mathcal{E}(n)$ is the corresponding eigenvalue.
This property is preserved under the affine transformation of $\eta(x)$,
$\eta(x)\to a\eta(x)+b$ ($a,b$: constants).
As will be shown shortly, the functional form of the sinusoidal
coordinate is characterised by the {\em closure relation\/}
\eqref{closurerel1}, \eqref{closurerel1t} up to a multiplicative
and an additive constant.
We choose the additive constant to achieve
\begin{equation}
  \eta(0)=0.
  \label{etazero}
\end{equation}
The lower triangularity \eqref{lowtri} implies that the eigenvectors
can be obtained as polynomials in $\eta(x)$ up to normalisation:
\begin{equation}
  \begin{array}{l}
   \widetilde{\mathcal{H}}P_n(\eta(x))=\mathcal{E}(n)P_n(\eta(x)),\\[8pt]
   \mathcal{H}\phi_n(x)
   =\mathcal{E}(n)\phi_n(x),\quad \phi_n(x)=\phi_0(x)P_n(\eta(x)),
  \end{array}
  \quad n=0,1,2,\ldots, (N)\ ,\ldots.
\end{equation}
We choose the simple normalisation of $P_n(\eta(x))$
\begin{equation}
  P_n(0)=1,\quad n=0,1,2,\ldots, (N)\ ,\ldots,
  \label{Pzero}
\end{equation}
which is consistent with \eqref{psizero} and \eqref{etazero}.
Note that this also means
\begin{equation}
  P_0=1.
  \label{Pzero1}
\end{equation}
Thus the eigenvectors are completely specified and the eigenvalue
problem \eqref{eigenpr} is completely solved.
The orthogonality of the eigenvectors of the Jacobi matrix
$\mathcal{H}$ \eqref{genham}, \eqref{Jacobiform} implies
\begin{equation}
  \sum_x\phi_0(x)^2P_n(\eta(x))P_m(\eta(x))=\frac{1}{d_n^2}\,\delta_{n,m},
  \label{phi0^2PnPn=dn^2}
\end{equation}
in which the normalisation constants $\{d_n>0\}$ are to be calculated.
In terms of the normalised eigenvectors
\begin{equation}
  \hat{\phi}_0(x)=d_0\phi_0(x),\quad
  \hat{\phi}_n(x)=d_n\phi_0(x)P_n(\eta(x))
  =\hat{\phi}_0(x)\frac{d_n}{d_0}P_n(\eta(x)),
\end{equation}
the orthonormality relation reads
\begin{equation}
  \sum_xd_n\phi_0(x)P_n(\eta(x))\cdot d_m\phi_0(x)P_m(\eta(x))
  =\delta_{n,m}.
\end{equation}
This in turn implies the completeness relation
\begin{equation}
  \sum_nd_n\phi_0(x)P_n(\eta(x))\cdot d_n\phi_0(y)P_n(\eta(y))
  =\delta_{x,y},
\end{equation}
which states the simplest duality that the row eigenvectors are also 
orthogonal.
When written slightly differently
\begin{equation}
  \sum_nd_n^2P_n(\eta(x))P_n(\eta(y))
  =\frac{1}{\phi_0(x)^2}\,\delta_{x,y},
  \label{sumdnPnPn}
\end{equation}
it displays the duality with \eqref{phi0^2PnPn=dn^2} explicitly.

We have established that $P_n(\eta)$ is an orthogonal polynomial of
a discrete variable with the discrete measure $\{\phi_0(x)^2\}$ given
explicitly as \eqref{phi0=prodB/D}, \eqref{phizero2}.
They satisfy the three term recurrence relation
\begin{equation}
  \eta P_n(\eta)=A_nP_{n+1}(\eta)+B_nP_n(\eta)+C_nP_{n-1}(\eta)
  \label{Pthreeterm}
\end{equation}
with coefficients $A_n$, $B_n$ and $C_n$ ($C_0=0$).
The boundary condition $P_n(0)=1$ \eqref{Pzero} implies
\begin{equation}
  B_n=-(A_n+C_n).
\end{equation}
In section \ref{shapeclosuresection} we will derive the general
formulas of the coefficients $A_n$ and $C_n$ \eqref{Anform},
\eqref{Cnform} from the input functions $B(x)$ and $D(x)$ with
the help of various quantities implied by the closure relation.
This will provide complete specification of the polynomial $P_n(\eta)$.

In the next section, we will show that the explicit form of the
normalisation constants $\{d_n\}$ \eqref{phi0^2PnPn=dn^2} can be
easily obtained as the components of the ground state wavefunction
of the {\em dual\/} polynomial as seen clearly from
\eqref{phi0^2PnPn=dn^2} and \eqref{sumdnPnPn}.

\section{Dual Polynomials}
\label{dualsection}
\setcounter{equation}{0}

The {\em dual\/} polynomial is an important concept in the theory
of orthogonal polynomials of a discrete variable
\cite{nikiforov,askey,gasper,ismail,koeswart}.
Here we show that the dual polynomial arises naturally as the
solution of the original eigenvalue problem \eqref{eigenpr}
or \eqref{vequation} obtained in a different way.
Let us rewrite the similarity transformed eigenvalue problem
$\widetilde{\mathcal H}\,v(x)=\mathcal{E}\,v(x)$ \eqref{vequation}
into an explicit matrix form with the change of the notation
$v(x)\to{}^t(Q_0,Q_1,\ldots,Q_x,\ldots)$
\begin{equation}
  \sum_{y}\widetilde{\mathcal H}_{x,y}Q_y=\mathcal{E}Q_x,
  \quad x=0,1,\ldots,(N),\ldots.
\end{equation}
Because of the tridiagonality of $\widetilde{\mathcal H}$, it is
in fact a three term recurrence relation for $Q_x$ as a polynomial
in $\mathcal{E}$:
\begin{align}
  \mathcal{E}Q_x(\mathcal{E})
  =B(x)\bigl(Q_x(\mathcal{E})-Q_{x+1}(\mathcal{E})\bigr)
  &+D(x)\bigl(Q_x(\mathcal{E})-Q_{x-1}(\mathcal{E})\bigr),\n
  &x=0,1,\ldots,(N),\ldots .
  \label{dual3term}
\end{align}
Starting with the boundary (initial) condition
\begin{equation}
  Q_0=1,
  \label{Qzero}
\end{equation}
which is consistent with \eqref{psizero}, we determine
$Q_x(\mathcal{E})$ as a degree $x$ polynomial in $\mathcal{E}$.
It is easy to see
\begin{equation}
  Q_x(0)=1,\quad x=0,1,\ldots,(N),\ldots .
  \label{Qxzero}
\end{equation}
When $\mathcal{E}$ is replaced by the actual value of the $n$-th
eigenvalue $\mathcal{E}(n)$ \eqref{lowtri} in $Q_x(\mathcal{E})$,
we obtain the explicit form of the eigenvector
\begin{equation}
  \sum_{y}\widetilde{\mathcal{H}}_{x,y}Q_y(\mathcal{E}(n))
  =\mathcal{E}(n)Q_x(\mathcal{E}(n)),
  \quad x=0,1,\ldots,(N),\ldots.
\end{equation}
In the finite dimensional case (I),
$\{Q_0(\mathcal{E}),\ldots,Q_N(\mathcal{E})\}$ are determined by
\eqref{dual3term} for $x=0,\ldots,N-1$. The last equation
\begin{align}
  \mathcal{E}Q_N(\mathcal{E})=
  D(N)\bigl(Q_N(\mathcal{E})-Q_{N-1}(\mathcal{E})\bigr),
\end{align}
is the degree $N+1$ algebraic equation (characteristic equation)
for the determination of all the eigenvalues $\{\mathcal{E}(n)\}$.
This is a familiar situation encountered in quasi-exactly solvable
quantum mechanics \cite{Ush,turb,benddun,sasaki07}.
In the next section we will provide further two independent algebraic
methods for the determination of the eigenvalues $\{\mathcal{E}(n)\}$
based on the shape-invariance \cite{os456} and the exact Heisenberg
operator solution \cite{os7}, both of which are applicable to the
finite (I) as well as the infinite dimensional case (II).

We now have two expressions (polynomials) for the eigenvectors
of the problem \eqref{eigenpr} belonging to the eigenvalue
$\mathcal{E}(n)$; $P_n(\eta(x))$ and $Q_x(\mathcal{E}(n))$.
Due to the simplicity of the spectrum of the Jacobi matrix, they must
be equal up to a multiplicative factor $\alpha_n$,
\begin{equation}
  P_n(\eta(x))=\alpha_n Q_x(\mathcal{E}(n)),\quad
  x=0,1,\ldots,(N),\ldots,
\end{equation}
which turns out to be unity because of the boundary (initial)
condition at $x=0$ \eqref{etazero}, \eqref{Pzero}, \eqref{Qzero}
(or at $n=0$ \eqref{Aphi0=0}, \eqref{Pzero1}, \eqref{Qxzero});
\begin{align}
  P_n(\eta(0))=P_n(0)&=1=Q_0(\mathcal{E}(n)),\qquad\qquad\ \ \,
  n=0,1,\ldots,(N),\ldots,\\
  P_0(\eta(x))&=1=Q_x(0)=Q_x(\mathcal{E}(0)),\quad
  x=0,1,\ldots,(N),\ldots .
\end{align}

We have established that two polynomials, $\{P_n(\eta)\}$ and
its {\em dual\/} polynomial $\{Q_x(\mathcal{E})\}$, coincides at
the integer lattice points:
\begin{equation}
  P_n(\eta(x))=Q_x(\mathcal{E}(n)),\quad n=0,1,\ldots,(N),\ldots,
  \quad x=0,1,\ldots,(N),\ldots.
  \label{Duality}
\end{equation}
The completeness relation \eqref{sumdnPnPn}
\begin{align}
  \sum_nd_n^2P_n(\eta(x))P_n(\eta(y))
  =\sum_nd_n^2Q_x(\mathcal{E}(n))Q_y(\mathcal{E}(n))
  =\frac{1}{\phi_0(x)^2}\,\delta_{x,y},
  \label{Qnortho}
\end{align}
is now understood as the orthogonality relation of the dual
polynomial $Q_x(\mathcal{E})$, and the previous normalisation
constant $d_n^2$ is now the orthogonality measure.

The real symmetric (hermitian) matrix $\mathcal{H}_{x,y}$
\eqref{Jacobiform} can be expressed in terms of the complete set
of the eigenvalues and the corresponding normalised eigenvectors
\begin{align}
  \mathcal{H}_{x,y}&=\sum_n\mathcal{E}(n)\hat{\phi}_n(x)\hat{\phi}_n(y),\\
  \hat{\phi}_n(x)&=d_n\phi_0(x)P_n(\eta(x))=d_n\phi_0(x)Q_x(\mathcal{E}(n)).
\end{align}
The very fact that it is tridiagonal can be easily verified by using
the difference equation for the polynomial $P_n(\eta(x))$ or the
three term recurrence relation for $Q_x(\mathcal{E}(n))$.

Here is the list of the dual correspondence:
\begin{align}
  x&\leftrightarrow n,
  \label{dualityxn}\\
  \eta(x)&\leftrightarrow \mathcal{E}(n),\qquad
  \eta(0)=0\leftrightarrow\mathcal{E}(0)=0,\\
  B(x)&\leftrightarrow -A_n,\\
  D(x)&\leftrightarrow -C_n,\\
  \frac{\phi_0(x)}{\phi_0(0)}&\leftrightarrow\frac{d_n}{d_0}.
  \label{dualityphi0dn}
\end{align}
In the last expression, we inserted $\phi_0(0)=1$ \eqref{iniconphi0}
for symmetry.
It should be remarked that $B(x)$ and $D(x)$ govern the difference
equation for the polynomial $P_n(\eta)$, the solution of which
requires the knowledge of the sinusoidal coordinate $\eta(x)$.
The same quantities $B(x)$ and $D(x)$ specify the three term
recurrence of the dual polynomial $Q_x(\mathcal{E})$ without the
knowledge of the spectrum $\mathcal{E}(n)$. It is required for them
to be the eigenvectors of the eigenvalue problem \eqref{eigenpr}.
Likewise, $A_n$ and $C_n$ in \eqref{Pthreeterm} specify the polynomial
$P_n(\eta)$ without the knowledge of the sinusoidal coordinate.
As for the dual polynomial $Q_x(\mathcal{E}(n))$, $A_n$ and $C_n$
provide the difference equation (in $n$), the solution of which
needs the explicit form of $\mathcal{E}(n)$. Let us stress that it is
the eigenvalue problem \eqref{eigenpr} with the specific Hamiltonian
\eqref{genham} that determines the polynomial $P_n(\eta)$ and its dual
$Q_x(\mathcal{E})$, the spectrum $\mathcal{E}(n)$, the sinusoidal
coordinate $\eta(x)$ and the orthogonality measures $\phi_0(x)^2$
and $d_n^2$.

\section{Shape Invariance and Closure Relation}
\label{shapeclosuresection}
\setcounter{equation}{0}

Like the Hamiltonians for the ($q$-)Askey-scheme of hypergeometric
orthogonal polynomials with continuous measure, all the Hamiltonians
for the orthogonal polynomials of a discrete variable discussed
in this paper are endowed with two types of symmetries;
shape-invariance \cite{genden,os456} and closure relation \cite{os7}.
The former leads to the determination of the entire spectrum
\eqref{spectrumform} and the corresponding eigenvectors
\eqref{eigvecform}. In other words, it guarantees the exact
solvability in the Schr\"odinger picture.
The closure relation \eqref{closurerel1}, on the other hand, is
essential for the exact Heisenberg operator solution \eqref{Heisensol},
which in turn gives the annihilation/creation operators
\eqref{a^{(pm)}}. They help to determine the coefficients $A_n$ and
$C_n$ \eqref{Anform}, \eqref{Cnform} in the three term recurrence
relation of the polynomial $P_n(\eta)$ as well as the algebraic
method \eqref{alphapmE} to determine the eigenvalues
$\{\mathcal{E}(n)\}$.
In \S\ref{detetaBD} we show that the requirement of the closure
relation determines the sinusoidal coordinates $\eta(x)$ and $B(x)$
and $D(x)$.
In Appendix A, the possible forms of the Hamiltonians are determined
\eqref{Bxformula2}--\eqref{Dxformula2}.

\subsection{Shape Invariance}
\label{shapesection}

Shape invariance dictates the parameter dependence of the Hamiltonian.
Let us denote the set of parameters symbolically by $\bm{\lambda}$:
\begin{equation}
  \mathcal{H}(x;\bm{\lambda})=\mathcal{A}(x;\bm{\lambda})^{\dagger}
  \mathcal{A}(x;\bm{\lambda}), \quad \text{or}\quad
  \mathcal{H}(\bm{\lambda})=\mathcal{A}(\bm{\lambda})^{\dagger}
  \mathcal{A}(\bm{\lambda}).
\end{equation}
Shape invariance simply means
\begin{equation}
  \mathcal{A}(\bm{\lambda})\mathcal{A}(\bm{\lambda})^{\dagger}
  =\kappa\mathcal{A}(\bm{\lambda}+\bm{\delta})^{\dagger}
  \mathcal{A}(\bm{\lambda}+\bm{\delta})
  +\mathcal{E}(1\,;\bm{\lambda}),
  \label{shapeinv1}
\end{equation}
in which $\kappa$ is a positive constant and
$\mathcal{E}(1\,;\bm{\lambda})$ is the eigenvalue of the first
excited state $\mathcal{E}(1)>0$ with the explicit parameter
dependence. That is, the original Hamiltonian
$\mathcal{H}(\bm{\lambda})$ and the {\em associated\/} Hamiltonian
$\mathcal{A}(\bm{\lambda})\mathcal{A}(\bm{\lambda})^{\dagger}$
in Crum's \cite{crum} sense (or in the factorisation method
\cite{infhul}, or in the so-called super-symmetric quantum mechanics
\cite{susyqm}) have the same {\em shape\/} with a multiplicative
factor $\kappa$.\footnote{
The continuous $q$-Hermite polynomial \cite{askey,koeswart}
has Hamiltonian \eqref{Hq1}
with $V(x)=1/\left((1-z^2)(1-qz^2)\right)$, $z=e^{ix}$, which has
no shiftable parameter. The above shape invariance relation gives a
simple realisation of the $q$-oscillator algebra
$\mathcal{A}\mathcal{A}^\dagger=
q^{-1}\mathcal{A}^\dagger\mathcal{A}+\mathcal{E}_1$,
$\mathcal{E}_1=q^{-1}-1$, \cite{os11}.
}
The necessary and sufficient condition for \eqref{shapeinv1} is
\begin{align}
  \sqrt{B(x+1\,;\bm{\lambda})D(x+1\,;\bm{\lambda})}
  &=\kappa
  \sqrt{B(x\,;\bm{\lambda}+\bm{\delta})D(x+1\,;\bm{\lambda}+\bm{\delta})},
  \label{shapeinv1cond1}\\
  B(x\,;\bm{\lambda})+D(x+1\,;\bm{\lambda})
  &=\kappa\bigl(
  B(x\,;\bm{\lambda}+\bm{\delta})+D(x\,;\bm{\lambda}+\bm{\delta})\bigr)
  +\mathcal{E}(1\,;\bm{\lambda}).
  \label{shapeinv1cond2}
\end{align}
For the finite dimensional case (I), the new (associated) Hamiltonian
$\mathcal{H}(\bm{\lambda}+\bm{\delta})
=\mathcal{A}(\bm{\lambda}+\bm{\delta})^{\dagger}
\mathcal{A}(\bm{\lambda}+\bm{\delta})$
is of $N$ dimensional, and \eqref{shapeinv1cond1} holds for
$0\leq x\leq N-1$.

It is rather straightforward to extract the necessary information,
as done explicitly for the ($q$-)Askey-scheme of hypergeometric
polynomials with continuous measure \cite{os456}.
The spectrum is simply generated by $\mathcal{E}(1,\bm{\lambda})$:
\begin{align}
  &\mathcal{E}(n\,;\bm{\lambda})
  =\sum_{s=0}^{n-1}\kappa^s\mathcal{E}(1\,;\bm{\lambda}+s\bm{\delta}),
  \label{spectrumform}
\end{align}
and the corresponding eigenvectors are generated from the known form
of the ground state eigenvector $\phi_0$ \eqref{phi0=prodB/D} together
with the multiple action of the successive $\mathcal{A}^\dagger$
operator:
\begin{align}
  &\phi_n(x\,;\bm{\lambda})\propto
  \mathcal{A}(\bm{\lambda})^{\dagger}
  \mathcal{A}(\bm{\lambda}+\bm{\delta})^{\dagger}
  \mathcal{A}(\bm{\lambda}+2\bm{\delta})^{\dagger}\cdots
  \mathcal{A}(\bm{\lambda}+(n-1)\bm{\delta})^{\dagger}
  \phi_0(x\,;\bm{\lambda}+n\bm{\delta}),
  \label{eigvecform}
\end{align}
which is related to a Rodrigues type formula. It should be stressed
that shape invariance, \eqref{spectrumform} guarantees the monotonous
increasing property (within the valid parameter range) of
$\mathcal{E}(n)$.
For all the polynomials discussed in this paper, the spectrum is
either linear or quadratic in $n$ or $q$-quadratic in $n$:
\begin{alignat}{2}
  \text{(\romannumeral1)}:&\quad&\mathcal{E}(n)&=n,
  \label{matheform1}\\
  \text{(\romannumeral2)}:&\quad&\mathcal{E}(n)&=\epsilon n(n+\alpha),
  \qquad \qquad \quad\
  \epsilon=\Bigl\{\begin{array}{ll}
  1&\text{for }\ \alpha>-1,\\
  -1&\text{for }\ \alpha<-N,
  \end{array}\\
  \text{(\romannumeral3)}:&\quad&\mathcal{E}(n)&=1-q^n,
  \label{matheform3}\\
  \text{(\romannumeral4)}:&\quad&\mathcal{E}(n)&=q^{-n}-1,
  \label{matheform4}\\
  \text{(\romannumeral5)}:&\quad&\mathcal{E}(n)&=\epsilon
  (q^{-n}-1)(1-\alpha q^n), \quad
  \epsilon=\Bigl\{\begin{array}{ll}
  1&\text{for }\ \alpha<q^{-1},\\
  -1&\text{for }\ \alpha>q^{-N}.
  \end{array}
  \label{matheform5}
\end{alignat}
As we will show shortly in the next subsection, the requirement of
closure relation restricts the structure of the sinusoidal coordinate
$\eta(x)$, which has the same form as above with $n$ replaced by $x$,
\eqref{etaform1}--\eqref{etaform5} reflecting the duality
$x\leftrightarrow n$, $\eta(x)\leftrightarrow\mathcal{E}(n)$.

\paragraph{Forward and Backward Shift Operators}

The counterparts of the $\mathcal{A}$ and $\mathcal{A}^\dagger$
operators in the polynomial space, called the forward and backward
shift operators, play an important role in the theory of orthogonal
polynomials \cite{koeswart,askey}.
Since the shift in $x$ is closely related with the parameter shift,
we introduce an auxiliary function
$\varphi(x)=\varphi(x\,;\bm{\lambda})>0$, $\varphi(0)=1$ to absorb
the effects of the parameter shift in the ground state wave function:
\begin{equation}
  \varphi(x\,;\bm{\lambda})\eqdef
  \sqrt{\frac{B(0\,;\bm{\lambda})}{B(x\,;\bm{\lambda})}}\,
  \frac{\phi_0(x\,;\bm{\lambda}+\bm{\delta})}{\phi_0(x\,;\bm{\lambda})}.
  \label{phi=Bvphiphi}
\end{equation}
In the finite dimensional case (I), $B(N\,;\bm{\lambda})=0$ but it is
canceled by the zero in $\phi_0(N\,;\bm{\lambda}+\bm{\delta})$ and
$\varphi$ is well defined. With the help of \eqref{phi0/phi0=B/D},
it can be rewritten as
\begin{equation}
  \varphi(x\,;\bm{\lambda})=
  \sqrt{\frac{B(0\,;\bm{\lambda})}{D(x+1\,;\bm{\lambda})}}\,
  \frac{\phi_0(x\,;\bm{\lambda}+\bm{\delta})}{\phi_0(x+1\,;\bm{\lambda})}.
  \label{phi=Dvphiphi}
\end{equation}
Explicitly it reads
\begin{equation}
  \varphi(x\,;\bm{\lambda})
  =\sqrt{\frac{B(0\,;\bm{\lambda})}{B(x\,;\bm{\lambda})}
  \prod_{y=0}^{x-1}
  \frac{B(y\,;\bm{\lambda}+\bm{\delta})}{B(y\,;\bm{\lambda})}
  \frac{D(y+1\,;\bm{\lambda})}{D(y+1\,;\bm{\lambda}+\bm{\delta})}}.
  \label{varphi}
\end{equation}
For the Hamiltonians \eqref{H1}, \eqref{Hq1} with a continuous variable
$x$, the corresponding $\varphi(x)$ is given in \cite{os7}.

Let us introduce the {\em forward shift operator\/}
$\mathcal{F}=\mathcal{F}(\bm{\lambda})$, and
the {\em backward shift operator\/}
$\mathcal{B}=\mathcal{B}(\bm{\lambda})$ as
\begin{align}
  \widetilde{\mathcal{H}}(\bm{\lambda})
  &=\mathcal{B}(\bm{\lambda})\mathcal{F}(\bm{\lambda}),\\
  \mathcal{F}(\bm{\lambda})&\eqdef
  \sqrt{B(0\,;\bm{\lambda})}\,
  \phi_0(x\,;\bm{\lambda}+\bm{\delta})^{-1}\circ
  \mathcal{A}(\bm{\lambda})\circ\phi_0(x\,;\bm{\lambda}),\\
  \mathcal{B}(\bm{\lambda})&\eqdef
  \frac{1}{\sqrt{B(0\,;\bm{\lambda})}}\,
  \phi_0(x\,;\bm{\lambda})^{-1}\circ
  \mathcal{A}(\bm{\lambda})^{\dagger}\circ
  \phi_0(x\,;\bm{\lambda}+\bm{\delta}).
\end{align}
With the help of \eqref{phi=Bvphiphi}--\eqref{phi=Dvphiphi}
we obtain  expressions
\begin{align}
  \mathcal{F}(\bm{\lambda})&=
  B(0\,;\bm{\lambda})\varphi(x\,;\bm{\lambda})^{-1}(1-e^{\partial}),\\
  \mathcal{B}(\bm{\lambda})&=
  \frac{1}{B(0\,;\bm{\lambda})}
  \bigl(B(x\,;\bm{\lambda})-D(x\,;\bm{\lambda})e^{-\partial}\bigr)
  \varphi(x\,;\bm{\lambda}).
\end{align}
With $B(0\,;\bm{\lambda})$ appearing in many places, these formulas
look rather contrived. However, their action on the polynomial
$P_n(\eta)=P_n(\eta(x\,;\bm{\lambda})\,;\bm{\lambda})$ is rather simple:
\begin{align}
  \mathcal{F}(\bm{\lambda})P_n(\eta(x\,;\bm{\lambda})\,;\bm{\lambda})
  &=f_n(\bm{\lambda})P_{n-1}(\eta(x\,;\bm{\lambda}+\bm{\delta})\,;
  \bm{\lambda}+\bm{\delta}),
  \label{FPn}\\
  \mathcal{B}(\bm{\lambda})P_n(\eta(x\,;\bm{\lambda}+\bm{\delta})\,;
  \bm{\lambda}+\bm{\delta})
  &=b_n(\bm{\lambda})P_{n+1}(\eta(x\,;\bm{\lambda})\,;\bm{\lambda}),
  \label{BPn}
\end{align}
in which $f_n(\bm{\lambda})$ and $b_n(\bm{\lambda})$ are constants
satisfying
$\mathcal{E}(n\,;\bm{\lambda})=b_{n-1}(\bm{\lambda})f_n(\bm{\lambda})$.
In fact, for all the examples given in \S\ref{examples}, they take
simple values,
\begin{equation}
  f_n(\bm{\lambda})=\mathcal{E}(n\,;\bm{\lambda}),\quad
  b_n(\bm{\lambda})=1.
\end{equation}
Moreover, the auxiliary function $\varphi(x)$ is related to the
sinusoidal coordinate $\eta(x)$ as
\begin{equation}
  \varphi(x)=\frac{\eta(x+1)-\eta(x)}{\eta(1)},
\end{equation}
which gives
\begin{equation}
  \frac{\eta(x)}{\eta(1)}=\sum_{y=0}^{x-1}\varphi(y).
\end{equation}
The action of $\mathcal{A}$, $\mathcal{A}^{\dagger}$ on the
eigenfunction $\phi_n(x)$  is
\begin{align}
   \mathcal{A}(\bm{\lambda})\phi_n(x\,;\bm{\lambda})
   &=\frac{1}{\sqrt{B(0\,;\bm{\lambda})}}\,
   f_n(\bm{\lambda})\phi_{n-1}(x\,;\bm{\lambda}+\bm{\delta}),
   \label{Aphin}\\
   \mathcal{A}(\bm{\lambda})^{\dagger}\phi_n(x\,;\bm{\lambda}+\bm{\delta})
   &=\sqrt{B(0\,;\bm{\lambda})}\,
   b_n(\bm{\lambda})\phi_{n+1}(x\,;\bm{\lambda}).
   \label{Adphin}
\end{align}

In terms of the forward and backward shift operators, the shape
invariance condition \eqref{shapeinv1} reads
\begin{equation}
  \mathcal{F}(\bm{\lambda})\mathcal{B}(\bm{\lambda})
  =\kappa\mathcal{B}(\bm{\lambda}+\bm{\delta})
  \mathcal{F}(\bm{\lambda}+\bm{\delta})
  +\mathcal{E}(1\,;\bm{\lambda})
  \label{shapeinv2}
\end{equation}
The necessary and sufficient condition for \eqref{shapeinv2} are
\begin{align}
  B(x+1\,;\bm{\lambda})\varphi(x+1\,;\bm{\lambda})
  &=\kappa B(x\,;\bm{\lambda}+\bm{\delta})
  \varphi(x\,;\bm{\lambda}),
  \label{shapeinvcond1}\\
  D(x\,;\bm{\lambda})\varphi(x-1\,;\bm{\lambda})
  &=\kappa D(x\,;\bm{\lambda}+\bm{\delta})
  \varphi(x\,;\bm{\lambda}),
  \label{shapeinvcond2}\\
  B(x\,;\bm{\lambda})+D(x+1\,;\bm{\lambda})
  &=\kappa\bigl(B(x\,;\bm{\lambda}+\bm{\delta})
  +D(x\,;\bm{\lambda}+\bm{\delta})\bigr)+\mathcal{E}(1\,;\bm{\lambda}),
  \label{shapeinvcond3}
\end{align}
which look simpler than \eqref{shapeinv1cond1} and \eqref{shapeinv1cond2}
thanks to the presence of $\varphi(x)$.

\subsection{Closure Relation}
\label{closuresection}

The other symmetry of the Hamiltonians for the orthogonal polynomials
is called the {\em closure relation\/}.
It is realised if the $\mathcal{H}=\mathcal{H}(\bm{\lambda})$ and
the {\em sinusoidal coordinate\/} $\eta=\eta(x\,;\bm{\lambda})$ satisfy
\begin{equation}
  [\mathcal{H},[\mathcal{H},\eta]\,]=\eta\,R_0(\mathcal{H})
  +[\mathcal{H},\eta]\,R_1(\mathcal{H})+R_{-1}(\mathcal{H}),
  \label{closurerel1}
\end{equation}
in which $R_i(z)=R_i(z\,;\bm{\lambda})$ is a polynomial in $z$.
The naming is due to the fact that
$\eta(x)+R_{-1}(\mathcal{H})/R_{0}(\mathcal{H})$ undergoes a sinusoidal
motion with frequency $\sqrt{R_{0}(\mathcal{H})}$ in the classical
dynamics limit \cite{os7}.
By similarity transformation in terms of $\phi_0(x\,;\bm{\lambda})$,
\eqref{closurerel1} reads
\begin{equation}
  [\widetilde{\mathcal{H}},[\widetilde{\mathcal{H}},\eta]\,]
  =\eta\,R_0(\widetilde{\mathcal{H}})
  +[\widetilde{\mathcal{H}},\eta]\,R_1(\widetilde{\mathcal{H}})
  +R_{-1}(\widetilde{\mathcal{H}}).
  \label{closurerel1t}
\end{equation}
The l.h.s. consists of the operators $e^{2\partial}$, $e^{\partial}$,
$1$, $e^{-\partial}$, $e^{-2\partial}$ and $\mathcal{H}$ contains
$e^{\pm\partial}$.
Thus $R_i$ can be parametrised as
\begin{equation}
  R_1(z)=r_1^{(1)}z+r_1^{(0)},\quad
  R_0(z)=r_0^{(2)}z^2+r_0^{(1)}z+r_0^{(0)},\quad
  R_{-1}(z)=r_{-1}^{(2)}z^2+r_{-1}^{(1)}z+r_{-1}^{(0)}.
  \label{Ricoeff}
\end{equation}
The coefficients $r_i^{(j)}$ depend on the overall normalisation of
the Hamiltonian and the sinusoidal coordinate:
$r_1^{(j)}\propto B(0)^{1-j}$, $r_0^{(j)}\propto B(0)^{2-j}$,
$r_{-1}^{(j)}\propto \eta(1)B(0)^{2-j}$.

\paragraph{Heisenberg operator solution and Annihilation and Creation
operators}

The closure relation \eqref{closurerel1} enables us to express any
multiple commutator
\[
  [\mathcal{H},[\mathcal{H},\cdots,[\mathcal{H},\eta(x)]
  \!\cdot\!\cdot\cdot]]
\]
as a linear combination of the operators $\eta(x)$ and
$[\mathcal{H},\eta(x)]$ with coefficients depending on the Hamiltonian
$\mathcal{H}$ only.
The exact Heisenberg operator solution for the sinusoidal coordinate
$\eta(x)$ \cite{os7} is given by ($t$ is the time variable):
\begin{align}
  &e^{it\mathcal{H}}\eta(x)e^{-it\mathcal{H}}
  =a^{(+)}e^{i\alpha_+(\mathcal{H})t}+a^{(-)}e^{i\alpha_-(\mathcal{H})t}
  -R_{-1}(\mathcal{H})R_0(\mathcal{H})^{-1},
  \label{Heisensol}\\
  &\alpha_{\pm}(z)\eqdef\tfrac12\bigl(R_1(z)
  \pm\sqrt{R_1(z)^2+4R_0(z)}\,\bigr),
  \label{alphapmexp}\\
  &\qquad\qquad
  R_1(z)=\alpha_+(z)+\alpha_-(z),\quad
  R_0(z)=-\alpha_+(z)\alpha_-(z),
  \label{R1R0}\\
  &a^{(\pm)}\eqdef\pm\Bigl([\mathcal{H},\eta(x)]-\bigl(\eta(x)
  +R_{-1}(\mathcal{H})R_0(\mathcal{H})^{-1}\bigr)\alpha_{\mp}(\mathcal{H})
  \Bigr)
  \bigl(\alpha_+(\mathcal{H})-\alpha_-(\mathcal{H})\bigr)^{-1}
  \label{a^{(pm)}}\\
  &\phantom{a^{(\pm)}}=
  \pm\bigl(\alpha_+(\mathcal{H})-\alpha_-(\mathcal{H})\bigr)^{-1}
  \Bigl([\mathcal{H},\eta(x)]+\alpha_{\pm}(\mathcal{H})\bigl(\eta(x)
  +R_{-1}(\mathcal{H})R_0(\mathcal{H})^{-1}\bigr)\Bigr).
\end{align}
The positive and negative frequency parts of the Heisenberg operator
solution, $a^{(-)}$ and $a^{(+)}$, are the {\em annihilation\/} and
{\em creation\/} operators, which are hermitian conjugate to each
other $a^{(+)\,\dagger}=a^{(-)}$. Applying  \eqref{Heisensol} to the
eigenvector $\phi_n(x)=\phi_0(x)P_n(\eta(x))$ and using the three
term recurrence relation for the polynomial $P_n(\eta)$
\eqref{Pthreeterm}, we obtain \cite{os7}
\begin{align}
  a^{(+)}\phi_n(x)&=A_n\phi_{n+1}(x),\\
  a^{(-)}\phi_n(x)&=C_n\phi_{n-1}(x).
\end{align}
The constant part in \eqref{Heisensol} is related to $B_n$:
\begin{equation}
  R_{-1}(\mathcal{E}(n))R_0(\mathcal{E}(n))^{-1}=-B_n=A_n+C_n.
  \label{Rm1E}
\end{equation}
The exact energy eigenvalues $\{\mathcal{E}(n)\}$ can be obtained
as the solution of the non-linear equation
\begin{equation}
  \alpha_{\pm}(\mathcal{E}(n))=\mathcal{E}(n\pm 1)-\mathcal{E}(n)
  \label{alphapmE}
\end{equation}
starting with $\mathcal{E}(0)=0$ \cite{os7}. The above relation
simplifies the form of the annihilation-creation operators applied
to $\phi_n(x)$
\begin{align}
  &a^{(\pm)}\phi_n(x)\\
  &=\frac{\pm1}{\mathcal{E}(n+1)-\mathcal{E}(n-1)}
  \Bigl([\mathcal{H},\eta(x)]
  +\bigl(\mathcal{E}(n)-\mathcal{E}(n\mp 1)\bigr)\eta(x)
  +\frac{R_{-1}(\mathcal{E}(n))}{\mathcal{E}(n\pm 1)-\mathcal{E}(n)}\Bigr)
  \phi_n(x).\nonumber
\end{align}
This is the eigenvector version of the raising/lowering operators or
the structure relation for orthogonal polynomials \cite{koorn}.
The polynomial version is obtained by the similarity transformation
in terms of $\phi_0(x)$,
$\tilde{a}^{(\pm)}\eqdef\phi_0(x)^{-1}\circ a^{(\pm)}\circ\phi_0(x)$,
\begin{align}
  &\tilde{a}^{(\pm)}P_n(\eta(x))\\
  &=\frac{\pm1}{\mathcal{E}(n+1)-\mathcal{E}(n-1)}
  \Bigl([\widetilde{\mathcal{H}},\eta(x)]
  +\bigl(\mathcal{E}(n)-\mathcal{E}(n\mp 1)\bigr)\eta(x)
  +\frac{R_{-1}(\mathcal{E}(n))}{\mathcal{E}(n\pm 1)-\mathcal{E}(n)}\Bigr)
  P_n(\eta(x)).\nonumber
\end{align}

\paragraph{Determination of $A_n$ and $C_n$}

The coefficients $A_n$ and $C_n$ of the three term recurrence
relation of the eigenpolynomial $P_n(\eta)$ can be determined
elementarily from the above expression for $B_n$ \eqref{Rm1E}
and the duality.
The starting point is the three term recurrence relations of the
polynomial $P_n(\eta)$ and its dual $Q_x(\mathcal{E})$:
\begin{align}
  \eta(x)P_n(\eta(x))&=A_n\bigl(P_{n+1}(\eta(x))-P_{n}(\eta(x))\bigr)
  +C_n\bigl(P_{n-1}(\eta(x))-P_{n}(\eta(x))\bigr),
  \label{sabuneq}\\
  \mathcal{E}(n)Q_x(\mathcal{E}(n))
  &=B(x)\bigl(Q_{x}(\mathcal{E}(n))-Q_{x+1}(\mathcal{E}(n))\bigr)
  +D(x)\bigl(Q_{x}(\mathcal{E}(n))-Q_{x-1}(\mathcal{E}(n))\bigr).
  \label{3term}
\end{align}
By putting $x=0$ in \eqref{3term} and using $Q_0=1$, $D(0)=0$, we obtain
\begin{equation}
  P_n(\eta(1))= Q_1(\mathcal{E}(n))=\frac{\mathcal{E}(n)-B(0)}{-B(0)},
  \label{mP1}
\end{equation}
in which the first equality is due to the duality \eqref{Duality}.
Next, by putting $x=1$ in \eqref{sabuneq} and using \eqref{mP1},
we obtain
\begin{equation}
  \eta(1)(\mathcal{E}(n)-B(0))=A_n(\mathcal{E}(n+1)-\mathcal{E}(n))+
  C_n(\mathcal{E}(n-1)-\mathcal{E}(n)).
  \label{ancneq}
\end{equation}
The two equations \eqref{Rm1E} and \eqref{ancneq} for $n\ge1$ give
\begin{align}
  A_n&=\frac{\frac{R_{-1}(\mathcal{E}(n))}{R_0(\mathcal{E}(n))}
  \bigl(\mathcal{E}(n)-\mathcal{E}(n-1)\bigr)
  +\eta(1)\bigl(\mathcal{E}(n)-B(0)\bigr)}
  {\mathcal{E}(n+1)-\mathcal{E}(n-1)}\,,
  \label{Anform}\\
  C_n&=\frac{\frac{R_{-1}(\mathcal{E}(n))}{R_0(\mathcal{E}(n))}
  \bigl(\mathcal{E}(n)-\mathcal{E}(n+1)\bigr)
  +\eta(1)\bigl(\mathcal{E}(n)-B(0)\bigr)}
  {\mathcal{E}(n-1)-\mathcal{E}(n+1)}\,,
  \label{Cnform}
\end{align}
and for $n=0$ we obtain from \eqref{Rm1E}
\begin{equation}
  A_0=R_{-1}(0)R_0(0)^{-1}=r^{(0)}_{-1}/r^{(0)}_{0},\quad C_0=0.
  \label{A0exp}
\end{equation}
They are slightly simplified ($n\ge1$) in terms of $\alpha_{\pm}$:
\begin{align}
  A_n&=\frac{R_{-1}(\mathcal{E}(n))
  +\eta(1)\bigl(\mathcal{E}(n)-B(0)\bigr)\alpha_+(\mathcal{E}(n))}
  {\alpha_+(\mathcal{E}(n))
  \bigl(\alpha_+(\mathcal{E}(n))-\alpha_-(\mathcal{E}(n))\bigr)}\,,
  \label{Anformula2}\\
  C_n&=\frac{R_{-1}(\mathcal{E}(n))
  +\eta(1)\bigl(\mathcal{E}(n)-B(0)\bigr)\alpha_-(\mathcal{E}(n))}
  {\alpha_-(\mathcal{E}(n))
  \bigl(\alpha_-(\mathcal{E}(n))-\alpha_+(\mathcal{E}(n))\bigr)}.
  \label{Cnformula2}
\end{align}
Note that \eqref{ancneq} for $n=0$ gives another important relation
\begin{equation}
  A_0\mathcal{E}(1)+B(0)\eta(1)=0.
  \label{ABrel}
\end{equation}
This relation ensures that the expression for $C_n$ \eqref{Cnformula2}
vanishes for $n=0$.
Written differently, the above relation \eqref{ABrel} means an equality
\begin{equation}
  \frac{B(0)}{\mathcal{E}(1)}=\frac{-A_0}{\eta(1)}
  \label{ABrel2}
\end{equation}
between the two intensive quantities. That is, they are independent of
the overall normalisation of the Hamiltonian $\mathcal{H}$ of the
polynomial system and of $\mathcal{H}^{\text{dual}}$ of its dual.

\paragraph{Dual Hamiltonian and Determination of $d_n/d_0$}

By the duality \eqref{Duality}, the three term recurrence relation
for $P_n(\eta)$ \eqref{sabuneq} translates to the difference equation
in $n$ for the dual polynomial $Q_x(\mathcal{E}(n))$;
\begin{equation}
  \eta(x)Q_x(\mathcal{E}(n))=
  A_n\bigl(Q_x(\mathcal{E}(n+1))-Q_x(\mathcal{E}(n))\bigr)
  +C_n\bigl(Q_x(\mathcal{E}(n-1))-Q_x(\mathcal{E}(n))\bigr).
\end{equation}
With abuse of notation, it can be written as  an eigenvalue equation
for the operator  $\widetilde{\mathcal H}^{\text{dual}}$ with
eigenvalue $\eta(x)$:
\begin{align}
  &\widetilde{\mathcal H}^{\text{dual}}\eqdef
  -A_n(1-e^{\partial_n})-C_n(1-e^{-\partial_n}),\\
  &\widetilde{\mathcal{H}}^{\text{dual}}Q_x(\mathcal{E}(n))
  =\eta(x)Q_x(\mathcal{E}(n)).
\end{align}
The corresponding Hamiltonian matrix is parametrised by
$n,m=0,1,\ldots,(N),\ldots$,
\begin{align}
  &\mathcal{H}^{\text{dual}}\eqdef
  -\sqrt{-A_n}\,e^{\partial_n}\sqrt{-C_n}
  -\sqrt{-C_n}\,e^{-\partial_n}\sqrt{-A_n}
  -(A_n+C_n),\ \ A_n<0,\ C_n<0,
  \label{dualgenham}\\
  &(\mathcal{H}^{\text{dual}})_{n,m}=
  -\sqrt{A_nC_{n+1}}\,\delta_{n+1,m}
  -\sqrt{A_{n-1}C_n}\,\delta_{n-1,m}
  -(A_n+C_n)\,\delta_{n,m}.
\end{align}
Following the same path that led to $\phi_0(x)$ \eqref{phi0=prodB/D},
the orthogonality measure for the polynomial $P_n(\eta(x))$,
we obtain the formula for the orthogonality measure $d_n$ for the
dual polynomial $Q_x(\mathcal{E}(n))$:
\begin{equation}
  \frac{d_n}{d_0}=\sqrt{\prod_{m=0}^{n-1}\frac{A_m}{C_{m+1}}}.
  \label{dn/d0}
\end{equation}
It should be stressed that this same $d_n$ is the normalisation
constant of the polynomial $P_n(\eta(x))$ \eqref{phi0^2PnPn=dn^2}.
With the above formula, one only needs to  evaluate  $d_0$ explicitly.

\subsection{Determination of $\eta(x)$ and $B(x)+D(x)$}
\label{detetaBD}

In the previous subsection, consequences of the closure relation are
explored. Here we will consider the closure relation as algebraic
constraints for unspecified Hamiltonian ({\em i.e.} $B(x)$ and $D(x)$)
and as yet undetermined sinusoidal coordinate $\eta(x)$.
Under very mild assumptions (on top of the general conditions
\eqref{BDcondition}), $\eta(0)=0$, $\eta(x)>\eta(y)$ for $x>y$,
$x,y\in\mathbb{Z}_+$, we determine the possible forms of $\eta(x)$
and $B(x)+D(x)$. The characterisation of all the systems satisfying
the closure relation ({\em i.e.} determination of $B(x)$ and $D(x)$
separately) will be relegated to Appendix A.

The necessary and sufficient condition for the closure relation in
the polynomial space \eqref{closurerel1t} is
\begin{align}
  &\eta(x+2)-2\eta(x+1)+\eta(x)=r_0^{(2)}\eta(x)+r_{-1}^{(2)}
  +r_1^{(1)}\bigl(\eta(x+1)-\eta(x)\bigr),
  \label{crcond1}\\
  &\eta(x-2)-2\eta(x-1)+\eta(x)=r_0^{(2)}\eta(x)+r_{-1}^{(2)}
  +r_1^{(1)}\bigl(\eta(x-1)-\eta(x)\bigr),
  \label{crcond1p}\\
  &\bigl(\eta(x+1)-\eta(x)\bigr)\bigl(B(x+1)+D(x+1)-B(x)-D(x)\bigr)\n
  &\quad=-\bigl(r_0^{(2)}\eta(x)+r_{-1}^{(2)}\bigr)
  \bigl(B(x+1)+D(x+1)+B(x)+D(x)\bigr)
  -\bigl(r_0^{(1)}\eta(x)+r_{-1}^{(1)}\bigr)\n
  &\quad\phantom{=}\;
  -\bigl(\eta(x+1)-\eta(x)\bigr)
  \Bigl(r_1^{(1)}\bigl(B(x+1)+D(x+1)\bigr)+r_1^{(0)}\Bigr),
  \label{crcond2}\\
  &\bigl(\eta(x-1)-\eta(x)\bigr)\bigl(B(x-1)+D(x-1)-B(x)-D(x)\bigr)\n
  &\quad=-\bigl(r_0^{(2)}\eta(x)+r_{-1}^{(2)}\bigr)
  \bigl(B(x-1)+D(x-1)+B(x)+D(x)\bigr)
  -\bigl(r_0^{(1)}\eta(x)+r_{-1}^{(1)}\bigr)\n
  &\quad\phantom{=}\;
  -\bigl(\eta(x-1)-\eta(x)\bigr)
  \Bigl(r_1^{(1)}\bigl(B(x-1)+D(x-1)\bigr)+r_1^{(0)}\Bigr),
  \label{crcond2p}\\
  &-2\bigl(\eta(x+1)-\eta(x)\bigr)B(x)D(x+1)
  -2\bigl(\eta(x-1)-\eta(x)\bigr)D(x)B(x-1)\n
  &\quad=\bigl(r_0^{(2)}\eta(x)+r_{-1}^{(2)}\bigr)
  \Bigl(\bigl(B(x)+D(x)\bigr)^2+B(x)D(x+1)+D(x)B(x-1)\Bigr)\n
  &\quad\phantom{=}\;
  +\bigl(r_0^{(1)}\eta(x)+r_{-1}^{(1)}\bigr)\bigl(B(x)+D(x)\bigr)
  +r_0^{(0)}\eta(x)+r_{-1}^{(0)}\n
  &\quad\phantom{=}\;
  +r_1^{(1)}\Bigl(\bigl(\eta(x+1)-\eta(x)\bigr)B(x)D(x+1)
  +\bigl(\eta(x-1)-\eta(x)\bigr)D(x)B(x-1)\Bigr).
  \label{crcond3}
\end{align}
The first two equations \eqref{crcond1} and \eqref{crcond1p},
two above and below the diagonal of \eqref{closurerel1t},
determine $\eta(x)$. Then the next two equations \eqref{crcond2}
and \eqref{crcond2p}, one above and below the diagonal, give
$B(x)+D(x)$ as a function of $\eta(x)$. In Appendix A, we show that
the last equation \eqref{crcond3}, the diagonal part, determines
$B(x)$ and $D(x)$ separately.

By subtracting \eqref{crcond1} from \eqref{crcond1p}
(with $x$ replaced by $x+2$), we obtain
\begin{equation}
  0=(r_1^{(1)}-r_0^{(2)})\bigl(\eta(x+2)-\eta(x)\bigr),
\end{equation}
which gives a general constraint on the parameters
$r_0^{(2)}=r_1^{(1)}$.
By subtracting \eqref{crcond2} from \eqref{crcond2p}
(with $x$ replaced by $x+1$) and using the above constraint, we obtain
\begin{equation}
  0=(2r_1^{(0)}-r_0^{(1)})\bigl(\eta(x+1)-\eta(x)\bigr).
\end{equation}
With the two constraints
\begin{equation}
  r_0^{(2)}=r_1^{(1)},\quad r_0^{(1)}=2r_1^{(0)},
  \label{cransr}
\end{equation}
the first equation \eqref{crcond1} is reduced to a three term
recurrence relation for $\eta(x)$:
\begin{equation}
  \eta(x+2)-(2+r_1^{(1)})\eta(x+1)+\eta(x)=r_{-1}^{(2)},
  \label{crcond1pp}
\end{equation}
which can be solved by the input $\eta(0)=0$ and $\eta(1)$.
The solutions are classified by the roots of the `characteristic
equation':
\[
  t^2-(2+r_1^{(1)})t+1=0.
\]
In order to have two positive roots (necessary for the condition
$\eta(x)>\eta(y)$ for $x>y$), $r_1^{(1)}\ge0$ is needed.
For the case $r_1^{(1)}=0$, we have a quadratic solution:
\begin{equation}
  \eta(x)=\frac{1}{2}r_{-1}^{(2)}x
  \Bigl(x-1+\frac{2\eta(1)}{r_{-1}^{(2)}}\Bigr) \quad
  (\text{for}\ \ r_{-1}^{(2)}\neq 0);\qquad
  \eta(x)=\eta(1)x \quad (\text{for}\ \ r_{-1}^{(2)}=0).
\end{equation}
For the rest $r_1^{(1)}>0$, we have a $q$-quadratic solution.
It should be noted that the overall multiplicative factor is irrelevant.
After choosing $\eta(1)$ and other parameters properly, we arrive at
\begin{alignat}{2}
  \text{(\romannumeral1)}:&\quad&\eta(x)&=x,
  \label{etaform1}\\
  \text{(\romannumeral2)}:&\quad&\eta(x)&=\epsilon'x(x+d),
  \qquad \qquad \quad \
  \epsilon'=\Bigl\{\begin{array}{ll}
  1&\text{for }\ d>-1,\\[4pt]
  -1&\text{for }\ d<-N,
  \label{etaform2}
  \end{array}\\
  \text{(\romannumeral3)}:&\quad&\eta(x)&=1-q^x,\\
  \text{(\romannumeral4)}:&\quad&\eta(x)&=q^{-x}-1,
  \label{etaform4}\\
  \text{(\romannumeral5)}:&\quad&\eta(x)&=\epsilon'(q^{-x}-1)(1-dq^x),\quad
  \epsilon'=\Bigl\{\begin{array}{ll}
  1&\text{for }\ d<q^{-1},\\[4pt]
  -1&\text{for }\ d>q^{-N}.
  \end{array}
  \label{etaform5}
\end{alignat}
Here as usual $q$ is $ 0<q<1 $.
Similar conclusion was reached in a different context by
Atakishiyev-Rahman-Suslov \cite{atakishi1}.
This result can be considered as a difference equation version of
Bochner's theorem \cite{bochner, askeywil, vinzhed}.
In fact, the original setting of Bochner's theorem \cite{bochner} is
rather limited and its generalisation within the framework of
ordinary quantum mechanics is given in [10], which contains various
sinusoidal coordinates $\eta(x)=x$, $x^2$, $\cos x$, $\sinh x$, $e^{-x}$, etc.

In these cases the parameters in \eqref{cransr} read
\begin{alignat}{2}
  \text{(\romannumeral1)}:&\quad&&
  r_1^{(1)}=0,\quad r_{-1}^{(2)}=0,\\
  \text{(\romannumeral2)}:&\quad&&
  r_1^{(1)}=0,\quad r_{-1}^{(2)}=2\epsilon',\\
  \text{(\romannumeral3)}:&\quad&&
  r_1^{(1)}=(q^{-\frac12}-q^{\frac12})^2,\quad
  r_{-1}^{(2)}=-(q^{-\frac12}-q^{\frac12})^2,\\
  \text{(\romannumeral4)}:&\quad&&
  r_1^{(1)}=(q^{-\frac12}-q^{\frac12})^2,\quad
  r_{-1}^{(2)}=(q^{-\frac12}-q^{\frac12})^2,\\
  \text{(\romannumeral5)}:&\quad&&
  r_1^{(1)}=(q^{-\frac12}-q^{\frac12})^2,\quad
  r_{-1}^{(2)}=(q^{-\frac12}-q^{\frac12})^2\epsilon'(1+d),
\end{alignat}
and $\eta(x\pm 1)-\eta(x)$ are given by
\begin{alignat}{2}
  \text{(\romannumeral1)}:&\quad&&\eta(x\pm 1)-\eta(x)=\pm 1,
  \label{etapm1}\\
  \text{(\romannumeral2)}:&\quad&&
  \eta(x\pm 1)-\eta(x)=\pm\epsilon'(2x\pm 1+d),\\
  \text{(\romannumeral3)}:&\quad&&
  \eta(x\pm 1)-\eta(x)=(1-q^{\pm 1})q^x,\\
  \text{(\romannumeral4)}:&\quad&&
  \eta(x\pm 1)-\eta(x)=(q^{\mp 1}-1)q^{-x},\\
  \text{(\romannumeral5)}:&\quad&&
  \eta(x\pm 1)-\eta(x)=\epsilon'(q^{\mp 1}-1)q^{-x}(1-dq^{2x\pm 1}).
  \label{etapm5}
\end{alignat}

Solving \eqref{crcond2} and \eqref{crcond2p} is somewhat involved.
With the constraints \eqref{cransr} incorporated, \eqref{crcond2} reads,
\begin{align}
  &\bigl(\eta(x+1)-\eta(x)\bigr)(a_{x+1}-a_x)\nonumber\\
  &\quad=-\bigl(r_1^{(1)}\eta(x)+r_{-1}^{(2)}\bigr)(a_{x+1}+a_x)
  -\bigl(2r_1^{(0)}\eta(x)+r_{-1}^{(1)}\bigr)
  -\bigl(\eta(x+1)-\eta(x)\bigr)(r_1^{(1)}a_{x+1}+r_1^{(0)}),\nonumber
\end{align}
which can be rewritten as
\begin{align}
  &\bigl(\eta(x+1)-\eta(x)+r_1^{(1)}\eta(x+1)+r_{-1}^{(2)}\bigr)a_{x+1}
  +\bigl(-\eta(x+1)+\eta(x)+r_1^{(1)}\eta(x)+r_{-1}^{(2)}\bigr)a_x\n
  &\quad
  =-\bigl(\eta(x+1)+\eta(x)\bigr)r_1^{(0)}-r_{-1}^{(1)}.
  \label{ninefurther}
\end{align}
Here we have introduced the abbreviation
\begin{equation}
  a_x\eqdef B(x)+D(x).
\end{equation}
With the help of \eqref{crcond1pp}, \eqref{ninefurther} can be
simplified to
\begin{equation}
  \bigl(\eta(x+2)-\eta(x+1)\bigr)a_{x+1}+\bigl(\eta(x-1)-\eta(x)\bigr)a_x
  =-\bigl(\eta(x+1)+\eta(x)\bigr)r_1^{(0)}-r_{-1}^{(1)}.
  \label{crcond2pp}
\end{equation}
By multiplying  $\eta(x)-\eta(x+1)$ to \eqref{crcond2pp}, we obtain
a two term recurrence relation for $a_x=B(x)+D(x)$:
\begin{align}
  &\bigl(\eta(x+2)-\eta(x+1)\bigr)\bigl(\eta(x)-\eta(x+1)\bigr)a_{x+1}
  -\bigl(\eta(x+1)-\eta(x)\bigr)\bigl(\eta(x-1)-\eta(x)\bigr)a_x\n
  &\qquad
  =r_1^{(0)}\bigl(\eta(x+1)^2-\eta(x)^2\bigr)
  +r_{-1}^{(1)}\bigl(\eta(x+1)-\eta(x)\bigr),
\end{align}
which can be solved with the boundary (initial) condition
\begin{equation}
  a_0=B(0).
\end{equation}
The result is that the diagonal component of $\mathcal{H}$ and
$\widetilde{\mathcal{H}}$, $a_x\eqdef B(x)+D(x)$, multiplied by
$\bigl(\eta(x+1)-\eta(x)\bigr)\bigl(\eta(x-1)-\eta(x)\bigr)$ is a
quadratic polynomial in $\eta(x)$:
\begin{equation}
  \bigl(\eta(x+1)-\eta(x)\bigr)\bigl(\eta(x-1)-\eta(x)\bigr)a_x
  =r_1^{(0)}\eta(x)^2+r_{-1}^{(1)}\eta(x)
  +\eta(1)\eta(-1)a_0,
  \label{ansa}
\end{equation}
in which we have introduced a constant $\eta(-1)$
\begin{equation}
  \eta(-1)\eqdef r_{-1}^{(2)}-\eta(1),
\end{equation}
such that \eqref{crcond1pp} is formally satisfied.
This is a very important characterisation of $B(x)$ and $D(x)$
as we will see in the concrete examples in section \ref{examples}.
It should be remarked that the factor
$\bigl(\eta(x+1)-\eta(x)\bigr)\bigl(\eta(x-1)-\eta(x)\bigr)$
is also a quadratic polynomial in $\eta(x)$,
\begin{equation}
  -\bigl(\eta(x+1)-\eta(x)\bigr)\bigl(\eta(x-1)-\eta(x)\bigr)=
  r_1^{(1)}\eta(x)^2+2r_{-1}^{(2)}\eta(x)-\eta(1)\eta(-1),
  \label{(eta-eta)(eta-eta)}
\end{equation}
for all the cases listed in \eqref{etaform1}--\eqref{etaform5},
as can be easily verified from \eqref{etapm1}--\eqref{etapm5}.
Therefore we obtain
\begin{equation}
  a_x=-\frac{r_1^{(0)}\eta(x)^2+r_{-1}^{(1)}\eta(x)+\eta(1)\eta(-1)B(0)}
  {r_1^{(1)}\eta(x)^2+2r_{-1}^{(2)}\eta(x)-\eta(1)\eta(-1)}.
  \label{ansa2}
\end{equation}

The analysis of the remaining equation \eqref{crcond3}, which is
the diagonal part of the closure relation \eqref{closurerel1t}
is relegated to Appendix A.

\subsection{Dual Closure Relation}
\label{dualclosuresection}

Since the sinusoidal coordinate plays the central role in the
theory of orthogonal polynomials, understanding its general
properties is important.
Let us now consider the `dual' closure relation by interchanging
$\eta\leftrightarrow\mathcal{H}$ $(\Leftrightarrow\mathcal{E})$
in the closure relation \eqref{closurerel1}:
\begin{equation}
  [\eta,[\eta,\mathcal{H}]\,]=\mathcal{H}\,R_0^{\text{dual}}(\eta)
  +[\eta,\mathcal{H}]\,R_1^{\text{dual}}(\eta)+R_{-1}^{\text{dual}}(\eta),
  \label{closurerel2}
\end{equation}
in which $R_i^{\text{dual}}(\eta)$ are polynomials in $\eta$.
By similarity transformation in terms of $\phi_0(x)$, we obtain
\begin{equation}
  [\eta,[\eta,\widetilde{\mathcal{H}}]\,]=
  \widetilde{\mathcal{H}}\,R_0^{\text{dual}}(\eta)
  +[\eta,\widetilde{\mathcal{H}}]\,R_1^{\text{dual}}(\eta)
  +R_{-1}^{\text{dual}}(\eta),
  \label{closurerel2t}
\end{equation}
which is again a relationship among tridiagonal and diagonal matrices.
It is equivalent to the following set of three equations:
\begin{align}
  &\bigl(\eta(x+1)-\eta(x)\bigr)^2
  =R_0^{\text{dual}}(\eta(x+1))-
  \bigl(\eta(x+1)-\eta(x)\bigr)R_1^{\text{dual}}(\eta(x+1)),
  \label{cr2cond1}\\
  &\bigl(\eta(x-1)-\eta(x)\bigr)^2
  =R_0^{\text{dual}}(\eta(x-1))-\bigl(\eta(x-1)
  -\eta(x)\bigr)R_1^{\text{dual}}(\eta(x-1)),
  \label{cr2cond1p}\\
  &0=\bigl(B(x)+D(x)\bigr)R_0^{\text{dual}}(\eta(x))
  +R_{-1}^{\text{dual}}(\eta(x)).
  \label{cr2cond2}
\end{align}
These imply
\begin{align}
  R_1^{\text{dual}}(\eta(x))&=
  \bigl(\eta(x+1)-\eta(x)\bigr)+\bigl(\eta(x-1)-\eta(x)\bigr),
  \label{dualclcon1}\\
  R_0^{\text{dual}}(\eta(x))&=
  -\bigl(\eta(x+1)-\eta(x)\bigr)\bigl(\eta(x-1)-\eta(x)\bigr),\\
  R_{-1}^{\text{dual}}(\eta(x))&
  =-\bigl(B(x)+D(x)\bigr)R_0^{\text{dual}}(\eta(x)).
  \label{dualclcon3}
\end{align}
These expressions are the dual formulas of $\mathcal{E}(n)$,
\eqref{R1R0}, \eqref{Rm1E} and \eqref{alphapmE}.
The right hand sides of \eqref{dualclcon1}--\eqref{dualclcon3} are
in fact polynomials in $\eta(x)$:
\begin{align}
  R_1^{\text{dual}}(z)&=
  r_1^{(1)}z+r_{-1}^{(2)},
  \label{dualclR1}\\
  R_0^{\text{dual}}(z)&=
  r_1^{(1)}z^2+2r_{-1}^{(2)}z-\eta(1)\eta(-1),
  \label{dualclR0}\\
  R_{-1}^{\text{dual}}(z)&=
  r_1^{(0)}z^2+r_{-1}^{(1)}z+\eta(1)\eta(-1)B(0).
  \label{dualclR-1}
\end{align}
The first two relations are consequences of \eqref{crcond1pp}
(see also \eqref{(eta-eta)(eta-eta)}), and the last relation is
nothing but the main result \eqref{ansa} of subsection \ref{detetaBD}.
The coefficients of $R^{\text{dual}}_i(z)$ defined like as
\eqref{Ricoeff} satisfy the same constraints as \eqref{cransr}.
The number of independent parameters in $R^{\text{dual}}_i(z)$
is six, just the same as in $R_i(z)$ with \eqref{cransr}.
Among them, as mentioned after \eqref{Ricoeff}, the overall
normalisation of $\mathcal{H}$ and $\eta(x)$ are immaterial.
This reduces the number of essential parameters to four,
which corresponds to the most generic cases of the Racah
\S\ref{[KS1.2]} and the $q$-Racah polynomials \S\ref{[KS3.2]}.

\section{Concrete Examples}
\label{examples}
\setcounter{equation}{0}

Here we will discuss the properties of various orthogonal
polynomials in our scheme.
We stick to the standard notation as far as possible;
the inevitable deviation is caused by symmetry consideration
and our universal normalisation of polynomials \eqref{Pzero},
\eqref{Qxzero}.
When we use non-standard parametrisation (Racah, $q$-Racah, Hahn,
$q$-Hahn and their duals) we adopt different symbols so that no
confusion will ensue.
We are well aware of the shortcomings of deviating from the
standard notation; for example, the connection to other existing
concepts/quantities are blurred. However, we strongly believe
that presenting the entire theory of orthogonal polynomials of
a discrete variable in the hermitian matrix (or operator) form
has enough merits of elucidating their common structure to
compensate the shortcomings.
The orthogonal polynomials whose orthogonality relation is written
in terms of a $q$-integral (Jackson integral),
{\em e.g.} the Big $q$-Jacobi polynomial, etc. will
need different treatment and will be discussed elsewhere.
We start with the finite dimensional case (I) and then move to
the infinite dimensional case (II).
In each case the Askey-scheme of hypergeometric orthogonal
polynomials will be followed by the $q$-scheme polynomials,
in the order presented in the review of Koekoek and Swarttouw
\cite{koeswart}, with the number {\em e.g.} [KS3.2] attached
indicating the subsection there.
The subsection is named after the polynomial to be discussed in
the subsection. We explain various quantities and concepts in
some detail for the first examples, the ($q$-) Racah polynomials.
For the rest, mostly formulas only are given except for those
needing some attention and/or remarks.
As mentioned earlier, a good part of the explicit formulas for 
specific polynomials is already known in various contexts.

The only input is the Hamiltonian $\mathcal{H}$ or the two positive
functions $B(x)$ and $D(x)$. In order to define them, the parameters
must be specified. So we start from the parameters and their shift
properties under shape invariance. The other quantities are defined
and derived as explained in previous sections.
Most basic definitions and notation are recapitulated in Appendix B
for self-containedness.

\subsection{Finite Dimensional Case (I)}

In this subsection $x$ takes a finite range of values
\[
  x=0,1,2,\ldots,N.
\]

\subsubsection{Racah [KS1.2] (self-dual  with different parameters)}
\label{[KS1.2]}

As well-known, the Racah polynomial is the most general
hypergeometric orthogonal polynomial of a discrete variable.
All the other (non-$q$) polynomials are obtained by restriction
or limiting procedure.
The set of parameters $\bm{\lambda}$ and their shifting unit
$\bm{\delta}$ and the multiplication factor $\kappa$ \eqref{shapeinv1}
are:
\begin{equation}
  \bm{\lambda}=(a,b,c,d),\quad\bm{\delta}=(1,1,1,1),\quad \kappa=1,
\end{equation}
which are different from the standard ones $(\alpha,\beta,\gamma,\delta)$.
Let us define $\tilde{d}$, $\epsilon$ and $\epsilon'$
\begin{align}
  &\tilde{d}\eqdef a+b+c-d-1,\\
  &\epsilon\eqdef \Bigl\{\begin{array}{ll}
  1&\text{for }\ \tilde{d}>-1,\\
  -1&\text{for }\ \tilde{d}<-N,
  \end{array}\quad
  \epsilon'\eqdef \Bigl\{\begin{array}{ll}
  1&\text{for }\ d>-1,\\
  -1&\text{for }\ d<-N,
  \end{array}
\end{align}
and $B(x)$ and $D(x)$ as
\begin{align}
  &B(x\,;\bm{\lambda})
  =-\epsilon\,\frac{(x+a)(x+b)(x+c)(x+d)}{(2x+d)(2x+1+d)},\\
  &D(x\,;\bm{\lambda})
  =-\epsilon\,\frac{(x+d-a)(x+d-b)(x+d-c)x}{(2x-1+d)(2x+d)}.
\end{align}
The finiteness condition $B(N)=0$ can be achieved by
\begin{equation}
  a=-N\quad\text{or}\quad b=-N\quad\text{or}\quad c=-N.
\end{equation}
Since $B(x)$ and $D(x)$ are symmetric in $a$, $b$ and $c$, we assume
without loss of generality
\begin{equation}
  c=-N,\quad a\ge b.
\end{equation}
Then we restrict our argument to the following parameter ranges
in which $B(x)$ and $D(x)$ are positive:
\begin{alignat}{2}
  &d>0,\ a>N+d,\ 0<b<1+d&\ \ \Rightarrow\ \ &
  (\epsilon,\epsilon')=(1,1),\\
  &d<-2N,\ a>0,\ N+d<b<1-N&\ \ \Rightarrow\ \ &
  (\epsilon,\epsilon')=(1,-1),\\
  &d>0,\ 0<a<1+d,\ b<1-N&\ \ \Rightarrow\ \ &
  (\epsilon,\epsilon')=(-1,1),\\
  &d<-2N,\ N+d<a<1-N,\ b<1+d&\ \ \Rightarrow\ \ &
  (\epsilon,\epsilon')=(-1,-1).
\end{alignat}
The energy eigenvalue and the sinusoidal coordinate are
\begin{equation}
  \mathcal{E}(n\,;\bm{\lambda})=\epsilon n(n+\tilde{d}),\quad
  \eta(x\,;\bm{\lambda})=\epsilon'x(x+d).
\end{equation}
The polynomial is
\begin{align}
  &P_n(\eta(x\,;\bm{\lambda})\,;\bm{\lambda})
  ={}_4F_3\Bigl(
  \genfrac{}{}{0pt}{}{-n,\,n+\tilde{d},\,-x,\,x+d}
  {a,\,b,\,c}\Bigm|1\Bigr)\n
  &\phantom{P_n(\eta(x\,;\bm{\lambda})\,;\bm{\lambda})}
  =R_n(\epsilon'\eta(x\,;\bm{\lambda})\,;a-1,\tilde{d}-a,c-1,d-c),
\end{align}
in which ${}_4F_3$ is the standard hypergeometric series
\eqref{defhypergeom} and $R_n$ is the standard notation for the
Racah polynomial \cite{koeswart}. The orthogonality measure (or
the squared ground state wavefunction) $\phi_0(x)^2$ and the
normalisation constants $d_n^2$ \eqref{phi0^2PnPn=dn^2} are
obtained from \eqref{phi0=prodB/D} and \eqref{dn/d0} (with $A_n$
and $C_n$ given below),
\begin{align}
  &\phi_0(x\,;\bm{\lambda})^2=\frac{(a,b,c,d)_x}{(1+d-a,1+d-b,1+d-c,1)_x}\,
  \frac{2x+d}{d},\\
  &d_n(\bm{\lambda})^2=\frac{(a,b,c,\tilde{d})_n}
  {(1+\tilde{d}-a,1+\tilde{d}-b,1+\tilde{d}-c,1)_n}\,
  \frac{2n+\tilde{d}}{\tilde{d}}\n
  &\qquad\qquad\quad\times
  \frac{(-1)^N(1+d-a,1+d-b,1+d-c)_N}{(\tilde{d}+1)_N(d+1)_{2N}},
\end{align}
in which $(a)_n$ is the Pochhammer symbol \eqref{defPoch}.
The format for $d_n(\bm{\lambda})^2$, throughout this section,
consists of two parts separated by a $\times$ symbol:
$d_n^2=(d_n^2/d_0^2)\times d_0^2$ reflecting the duality
\eqref{dualityphi0dn}. The second part $d_0^2$ satisfies the
relation $\sum_x\phi_0(x)^2=1/d_0^2$.
The coefficients of the three term recurrence \eqref{Pthreeterm}
$A_n$ and $C_n$ for the polynomial $P_n(\eta)$ obtained from the
closure relation \eqref{Anform}--\eqref{Cnform} are:
\begin{align}
  &A_n(\bm{\lambda})=\epsilon'\,\frac{(n+a)(n+b)(n+c)(n+\tilde{d})}
  {(2n+\tilde{d})(2n+1+\tilde{d})},\\
  &C_n(\bm{\lambda})=\epsilon'\,
  \frac{(n+\tilde{d}-a)(n+\tilde{d}-b)(n+\tilde{d}-c)n}
  {(2n-1+\tilde{d})(2n+\tilde{d})}.
\end{align}
The three functions $R_1$, $R_0$ and $R_{-1}$ appearing in the
closure relation \eqref{closurerel1} are:
\begin{align}
  &R_1(z\,;\bm{\lambda})=2\epsilon,\\
  &R_0(z\,;\bm{\lambda})=4\epsilon z+\tilde{d}^{\,2}-1,\\
  &R_{-1}(z\,;\bm{\lambda})=2\epsilon' z^2
  +\epsilon\epsilon'\bigl(2(ab+bc+ca)-(1+d)(1+\tilde{d})\bigr)z
  +\epsilon'abc(\tilde{d}-1).
\end{align}
Under the parameter shift $\bm{\lambda}\to \bm{\lambda}+\bm{\delta}$
the functions $B(x)$ and $D(x)$ behave:
\begin{equation}
  B(x\,;\bm{\lambda}+\bm{\delta})=B(x+1\,;\bm{\lambda})\,
  \frac{2x+3+d}{2x+1+d}\,,\quad
  D(x\,;\bm{\lambda}+\bm{\delta})=D(x\,;\bm{\lambda})\,
  \frac{2x-1+d}{2x+1+d}.
\end{equation}
{}From this and \eqref{varphi}, $\varphi(x)$ is
\begin{equation}
  \varphi(x\,;\bm{\lambda})=\frac{2x+d+1}{d+1}\,.
\end{equation}
The constants $f_n(\bm{\lambda})$ and $b_n(\bm{\lambda})$ appearing
in \eqref{FPn}-\eqref{BPn} are
\begin{equation}
  f_n(\bm{\lambda})=\mathcal{E}(n\,;\bm{\lambda}),\quad
  b_n(\bm{\lambda})=1.
\end{equation}
By comparing $B(x)$, $D(x)$ with $A_n$ and $C_n$, $\phi_0(x)^2$
with $d_n^2$, it is clear that the theory is self-dual with the
parameter correspondence
$(a,b,c,d,\epsilon,\epsilon')\leftrightarrow
(a,b,c,\tilde{d},\epsilon',\epsilon)$.

\subsubsection{Hahn [KS1.5] (dual to dual Hahn \S\ref{[KS1.6]})}
\label{[KS1.5]}

The parameters $(a,b)$ are very slightly different from the standard
ones $(\alpha,\beta)$ ($\alpha+1=a, \beta+1=b$) for the Hahn
polynomial. For obvious reasons, we adopt the same parameters for
the Hahn and dual Hahn polynomials:
\begin{align}
  &\bm{\lambda}=(a,b,N),\quad
  \bm{\delta}=(1,1,-1),\quad \kappa=1,\\
  &\epsilon=\Bigl\{\begin{array}{ll}
  1&\text{for }\ a,b>0,\\
  -1&\text{for }\ a,b<1-N,
  \end{array}\\
  &B(x\,;\bm{\lambda})=\epsilon(x+a)(N-x),\quad
  D(x\,;\bm{\lambda})=\epsilon x(b+N-x),
  \label{hahnBD}\\
  &\mathcal{E}(n\,;\bm{\lambda})=\epsilon n(n+a+b-1),\quad
  \eta(x\,;\bm{\lambda})=x,
  \label{hahneeta}\\
  &P_n(\eta(x\,;\bm{\lambda})\,;\bm{\lambda})
  ={}_3F_2\Bigl(
  \genfrac{}{}{0pt}{}{-n,\,n+a+b-1,\,-x}
  {a,\,-N}\Bigm|1\Bigr)
  =Q_n(\eta(x\,;\bm{\lambda})\,;a-1,b-1,N),
  \label{hahnpolydef}\\
  &\phi_0(x\,;\bm{\lambda})^2
  =\frac{N!}{x!\,(N-x)!}\,\frac{(a)_x\,(b)_{N-x}}{(b)_N}\,,
  \label{hahnphi0}\\
  &d_n(\bm{\lambda})^2
  =\frac{N!}{n!\,(N-n)!}\,
  \frac{(a)_n\,(2n+a+b-1)(a+b)_N}{(b)_n\,(n+a+b-1)_{N+1}}
  \times\frac{(b)_N}{(a+b)_N}\,,
  \label{hahndn}\\
  &A_n(\bm{\lambda})=-\frac{(n+a)(n+a+b-1)(N-n)}{(2n-1+a+b)(2n+a+b)}\,,
  \label{hahnAC1}\\
  &C_n(\bm{\lambda})=-\frac{n(n+b-1)(n+a+b+N-1)}{(2n-2+a+b)(2n-1+a+b)}\,,
  \label{hahnAC2}\\
  &R_1(z\,;\bm{\lambda})=2\epsilon,\\
  &R_0(z\,;\bm{\lambda})=4\epsilon z+(a+b-2)(a+b),\\
  &R_{-1}(z\,;\bm{\lambda})=-\epsilon(2N-a+b)z-a(a+b-2)N,\\
  &B(x\,;\bm{\lambda}+\bm{\delta})=B(x+1\,;\bm{\lambda}),\quad
  D(x\,;\bm{\lambda}+\bm{\delta})=D(x\,;\bm{\lambda}),\\
  &\varphi(x\,;\bm{\lambda})=1,\quad
  f_n(\bm{\lambda})=\mathcal{E}(n\,;\bm{\lambda}),\quad
  b_n(\bm{\lambda})=1.
\end{align}
Obviously the Hahn and dual Hahn polynomials are dual to each other.
The standard parameters for the latter are $(\gamma,\delta)$.

\subsubsection{dual Hahn [KS1.6] (dual to \S\ref{[KS1.5]})}
\label{[KS1.6]}

The parameters $(a,b)$ are very slightly different from the standard
ones $(\gamma,\delta)$ ($\gamma+1=a,\delta+1=b$) for the dual Hahn
polynomial:
\begin{align}
  &\bm{\lambda}=(a,b,N),\quad
  \bm{\delta}=(1,0,-1),\quad \kappa=1,\\
  &\epsilon=\Bigl\{\begin{array}{ll}
  1&\text{for }\ a,b>0,\\
  -1&\text{for }\ a,b<1-N,
  \end{array}\\
  &B(x\,;\bm{\lambda})=\frac{(x+a)(x+a+b-1)(N-x)}
  {(2x-1+a+b)(2x+a+b)},
  \label{dualhahnBD1}\\
  &D(x\,;\bm{\lambda})=\frac{x(x+b-1)(x+a+b+N-1)}
  {(2x-2+a+b)(2x-1+a+b)},
  \label{dualhahnBD2}\\
  &\mathcal{E}(n\,;\bm{\lambda})=n,\quad
  \eta(x\,;\bm{\lambda})=\epsilon x(x+a+b-1),
  \label{dualhahneeta}\\
  &P_n(\eta(x\,;\bm{\lambda})\,;\bm{\lambda})
  ={}_3F_2\Bigl(
  \genfrac{}{}{0pt}{}{-n,\,x+a+b-1,\,-x}
  {a,\,-N}\Bigm|1\Bigr)
  =R_n(\epsilon\eta(x\,;\bm{\lambda})\,;a-1,b-1,N),
  \label{dualhanpolydef}\\
  &\phi_0(x\,;\bm{\lambda})^2
  =\frac{N!}{x!\,(N-x)!}\,
  \frac{(a)_x\,(2x+a+b-1)(a+b)_N}{(b)_x\,(x+a+b-1)_{N+1}\,}\,,
  \label{dualhahnphi0}\\
  &d_n(\bm{\lambda})^2
  =\frac{N!}{n!\,(N-n)!}\,\frac{(a)_n\,(b)_{N-n}}{(b)_N}
  \times\frac{(b)_{N}}{(a+b)_N},
  \label{dualhahndn}\\
  &A_n(\bm{\lambda})=-\epsilon(n+a)(N-n),\quad
  C_n(\bm{\lambda})=-\epsilon n(b+N-n),
  \label{dualhahnAC}\\
  &R_1(z\,;\bm{\lambda})=0,\quad
  R_0(z\,;\bm{\lambda})=1,\quad
  R_{-1}(z\,;\bm{\lambda})=2\epsilon z^2-\epsilon(2N-a+b)z-\epsilon aN,\\
  &B(x\,;\bm{\lambda}+\bm{\delta})=B(x+1\,;\bm{\lambda})\,
  \frac{2x+2+a+b}{2x+a+b},\ \,
  D(x\,;\bm{\lambda}+\bm{\delta})=D(x\,;\bm{\lambda})\,
  \frac{2x-2+a+b}{2x+a+b},\\
  &\varphi(x\,;\bm{\lambda})=\frac{2x+a+b}{a+b},\quad
  f_n(\bm{\lambda})=\mathcal{E}(n\,;\bm{\lambda}),\quad
  b_n(\bm{\lambda})=1.
\end{align}
With the present parametrisation, the following duality
\eqref{hahnpolydef} $\leftrightarrow$ \eqref{dualhanpolydef},
\eqref{hahneeta} $\leftrightarrow$ \eqref{dualhahneeta},
\eqref{hahnBD} $\leftrightarrow$ \eqref{dualhahnAC},
\eqref{hahnAC1}--\eqref{hahnAC2} $\leftrightarrow$
\eqref{dualhahnBD1}--\eqref{dualhahnBD2},
\eqref{hahnphi0} $\leftrightarrow$ \eqref{dualhahndn},
\eqref{hahndn} $\leftrightarrow$ \eqref{dualhahnphi0},
is obvious.

\subsubsection{Krawtchouk [KS1.10] (self-dual)}
\label{[KS1.10]}

\begin{align}
  &\bm{\lambda}=(p,N),\quad
  \bm{\delta}=(0,-1),\quad \kappa=1,\\
  &0<p<1,\\
  &B(x\,;\bm{\lambda})=p(N-x),\quad
  D(x\,;\bm{\lambda})=(1-p)x,\\
  &\mathcal{E}(n\,;\bm{\lambda})=n,\quad
  \eta(x\,;\bm{\lambda})=x,\\
  &P_n(\eta(x\,;\bm{\lambda})\,;\bm{\lambda})
  ={}_2F_1\Bigl(
  \genfrac{}{}{0pt}{}{-n,\,-x}{-N}\Bigm|p^{-1}\Bigr)
  =K_n(\eta(x\,;\bm{\lambda})\,;p,N),\\
  &\phi_0(x\,;\bm{\lambda})=
  \frac{N!}{x!\,(N-x)!}\Bigl(\frac{p}{1-p}\Bigr)^x,\\
  &d_n(\bm{\lambda})^2
  =\frac{N!}{n!\,(N-n)!}\Bigl(\frac{p}{1-p}\Bigr)^n\times(1-p)^N,\\
  &A_n(\bm{\lambda})=-p(N-n),\quad
  C_n(\bm{\lambda})=-(1-p)n,\\
  &R_1(z\,;\bm{\lambda})=0,\quad
  R_0(z\,;\bm{\lambda})=1,\quad
  R_{-1}(z\,;\bm{\lambda})=(2p-1)z-pN,\\
  &B(x\,;\bm{\lambda}+\bm{\delta})=B(x+1\,;\bm{\lambda}),\quad
  D(x\,;\bm{\lambda}+\bm{\delta})=D(x\,;\bm{\lambda}),\\
  &\varphi(x\,;\bm{\lambda})=1,\quad
  f_n(\bm{\lambda})=\mathcal{E}(n\,;\bm{\lambda}),\quad
  b_n(\bm{\lambda})=1.
\end{align}

\subsubsection{$q$-Racah [KS3.2] (self-dual with different parameters)}
\label{[KS3.2]}

This is the first example of the $q$-scheme of the orthogonal
polynomials. Among them the $q$-Racah polynomial is the most
general. The set of parameters $\bm{\lambda}$ is different from
the standard one $(\alpha,\beta,\gamma,\delta)$ in the same manner
as for the Racah polynomial.
The shifting unit $\bm{\delta}$ and the multiplication factor $\kappa$
\eqref{shapeinv1} are:
\begin{equation}
  q^{\bm{\lambda}}=(a,b,c,d),\quad\bm{\delta}=(1,1,1,1),
  \quad \kappa=q^{-1}.
\end{equation}
Here $q^{\bm{\lambda}}$ stands for
$q^{(\lambda_1,\lambda_2,\ldots)}=(q^{\lambda_1},q^{\lambda_2},\ldots)$.
They are shifted multiplicatively. By definition the parameter $q$
is not shifted. Let us introduce $\tilde{d}$, $\epsilon$, $\epsilon'$
\begin{align}
  &\tilde{d}=abcd^{-1}q^{-1},\\
  &\epsilon=\Bigl\{\begin{array}{ll}
  1&\text{for }\ \tilde{d}<q^{-1},\\
  -1&\text{for }\ \tilde{d}>q^{-N},
  \end{array}\quad
  \epsilon'=\Bigl\{\begin{array}{ll}
  1&\text{for }\ d<q^{-1},\\
  -1&\text{for }\ d>q^{-N}.
  \end{array}
\end{align}
The functions $B(x)$ and $D(x)$ are
\begin{align}
  &B(x\,;\bm{\lambda})
  =-\epsilon\,\frac{(1-aq^x)(1-bq^x)(1-cq^x)(1-dq^x)}
  {(1-dq^{2x})(1-dq^{2x+1})}\,,\\
  &D(x\,;\bm{\lambda})
  =-\epsilon \tilde{d}\,
  \frac{(1-a^{-1}dq^x)(1-b^{-1}dq^x)(1-c^{-1}dq^x)(1-q^x)}
  {(1-dq^{2x-1})(1-dq^{2x})}.
\end{align}
The finiteness condition $B(N)=0$ can be realised
by one of the following choices:
\begin{equation}
  a=q^{-N}\quad\text{or}\quad b=q^{-N}\quad\text{or}\quad c=q^{-N}.
\end{equation}
Thanks to the symmetry in $a$, $b$ and $c$, we choose for simplicity
\begin{equation}
  c=q^{-N},\quad 0<a\leq b,\quad 0<d,
\end{equation}
and consider the following parameter ranges in which $B(x)$ and
$D(x)$ are positive:
\begin{alignat}{2}
  &0<d<1,\ 0<a<q^Nd,\ qd<b<1&\ \ \Rightarrow\ \ &
  (\epsilon,\epsilon')=(1,1),\\
  &d>q^{-2N},\ 0<a<1,\ q^{1-N}<b<q^Nd&\ \ \Rightarrow\ \ &
  (\epsilon,\epsilon')=(1,-1),\\
  &0<d<1,\ qd<a<1,\ b>q^{1-N}&\ \ \Rightarrow\ \ &
  (\epsilon,\epsilon')=(-1,1),\\
  &d>q^{-2N},\ q^{1-N}<a<q^Nd,\ b>qd&\ \ \Rightarrow\ \ &
  (\epsilon,\epsilon')=(-1,-1).
\end{alignat}
In case some of the parameters $a$, $b$ and $d$ are non-positive,
the situation is slightly complicated but can be treated in a
similar way.
The energy eigenvalue, the sinusoidal coordinate and the
$q$-polynomial are:
\begin{align}
  &\mathcal{E}(n\,;\bm{\lambda})=\epsilon(q^{-n}-1)(1-\tilde{d}q^n),\quad
  \eta(x\,;\bm{\lambda})=\epsilon'(q^{-x}-1)(1-dq^x),\\
  &P_n(\eta(x\,;\bm{\lambda})\,;\bm{\lambda})
  ={}_4\phi_3\Bigl(
  \genfrac{}{}{0pt}{}{q^{-n},\,\tilde{d}q^n,\,q^{-x},\,dq^x}
  {a,\,b,\,c}\Bigm|q\,;q\Bigr)\n
  &\phantom{P_n(\eta(x\,;\bm{\lambda})\,;\bm{\lambda})}
  =R_n(\epsilon'(1+d+\eta(x\,;\bm{\lambda}))\,;
  aq^{-1},\tilde{d}a^{-1},cq^{-1},dc^{-1}),
\end{align}
in which ${}_4\phi_3$ is the basic hypergeometric series
\eqref{defqhypergeom}.
The orthogonality measure $\phi_0(x)^2$ and the normalisation constants
$d_n^2$ \eqref{phi0^2PnPn=dn^2} are obtained from \eqref{phi0=prodB/D}
and \eqref{dn/d0} (with $A_n$ and $C_n$ given below),
\begin{align}
  &\phi_0(x\,;\bm{\lambda})^2=\frac{(a,b,c,d\,;q)_x}
  {(a^{-1}dq,b^{-1}dq,c^{-1}dq,q\,;q)_x\,\tilde{d}^x}\,
  \frac{1-dq^{2x}}{1-d}\,,\\
  &d_n(\bm{\lambda})^2
  =\frac{(a,b,c,\tilde{d}\,;q)_n}
  {(a^{-1}\tilde{d}q,b^{-1}\tilde{d}q,c^{-1}\tilde{d}q,q\,;q)_n\,d^n}\,
  \frac{1-\tilde{d}q^{2n}}{1-\tilde{d}}\n
  &\qquad\qquad\quad\times
  \frac{(-1)^N(a^{-1}dq,b^{-1}dq,c^{-1}dq\,;q)_N\,\tilde{d}^Nq^{\frac12N(N+1)}}
  {(\tilde{d}q\,;q)_N(dq\,;q)_{2N}}\,,
\end{align}
in which $(a;q)_n$ is the $q$-Pochhammer symbol \eqref{defqPoch}.
The coefficients of the three term recurrence \eqref{Pthreeterm}
$A_n$ and $C_n$ for the polynomial $P_n(\eta)$ obtained from the
closure relation \eqref{Anform}--\eqref{Cnform} are:
\begin{align}
  &A_n(\bm{\lambda})
  =\epsilon'\,\frac{(1-aq^n)(1-bq^n)(1-cq^n)(1-\tilde{d}q^n)}
  {(1-\tilde{d}q^{2n})(1-\tilde{d}q^{2n+1})}\,,\\
  &C_n(\bm{\lambda})
  =\epsilon'd\,
  \frac{(1-a^{-1}\tilde{d}q^n)(1-b^{-1}\tilde{d}q^n)(1-c^{-1}\tilde{d}q^n)
  (1-q^n)}
  {(1-\tilde{d}q^{2n-1})(1-\tilde{d}q^{2n})}.
\end{align}
The three functions $R_1$, $R_0$ and $R_{-1}$ appearing in the closure
relation \eqref{closurerel1} are:
\begin{align}
  &R_1(z\,;\bm{\lambda})=(q^{-\frac12}-q^{\frac12})^2z',\quad
  z'\eqdef z+\epsilon(1+\tilde{d}),\\
  &R_0(z\,;\bm{\lambda})=(q^{-\frac12}-q^{\frac12})^2
  \bigl(z^{\prime\,2}-(q^{-\frac12}+q^{\frac12})^2\tilde{d}\,\bigr),\\
  &R_{-1}(z\,;\bm{\lambda})=(q^{-\frac12}-q^{\frac12})^2
  \Bigl(\epsilon'(1+d)z^{\prime\,2}
  -\epsilon\epsilon'\bigl(a+b+c+d+\tilde{d}+(ab+bc+ca)q^{-1}\bigr)z'\n
  &\qquad\qquad\qquad\qquad\qquad
  +\epsilon'\bigl((1-a)(1-b)(1-c)(1-\tilde{d}q^{-1})\n
  &\qquad\qquad\qquad\qquad\qquad\qquad
  +(a+b+c-1-d\tilde{d}+(ab+bc+ca)q^{-1})\bigr)(1+\tilde{d})\Bigr).
\end{align}
For the shifted parameters the functions $B$ and $D$ are:
\begin{equation}
  B(x\,;\bm{\lambda}+\bm{\delta})=B(x+1\,;\bm{\lambda})\,
  \frac{1-dq^{2x+3}}{1-dq^{2x+1}}\,,\quad
  D(x\,;\bm{\lambda}+\bm{\delta})=q^2D(x\,;\bm{\lambda})\,
  \frac{1-dq^{2x-1}}{1-dq^{2x+1}}.
\end{equation}
{}From this and \eqref{varphi}, $\varphi(x)$ is
\begin{equation}
  \varphi(x\,;\bm{\lambda})=\frac{q^{-x}-dq^{x+1}}{1-dq}\,.
\end{equation}
The constants $f_n(\bm{\lambda})$ and $b_n(\bm{\lambda})$ appearing
in \eqref{FPn}-\eqref{BPn} are
\begin{equation}
  f_n(\bm{\lambda})=\mathcal{E}(n\,;\bm{\lambda}),\quad
  b_n(\bm{\lambda})=1.
\end{equation}
By comparing $B(x)$, $D(x)$ with $A_n$ and $C_n$, $\phi_0(x)^2$
with $d_n^2$, it is evident that the theory is self-dual with
the parameter correspondence
$(a,b,c,d,\epsilon,\epsilon')\leftrightarrow
(a,b,c,\tilde{d},\epsilon',\epsilon)$.

\subsubsection{$q$-Hahn [KS3.6] (dual to \S\ref{[KS3.7]})}
\label{[KS3.6]}

Here and in the next subsection (dual $q$-Hahn) we will consider,
for simplicity, the positive parameter range of $a$ and $b$,
which are very slightly different from the standard ones
$(\alpha,\beta)$ ($\alpha q=a,\beta q=b$) for the Hahn and dual
Hahn polynomials.
The standard parameters for the latter are $(\gamma,\delta)$.
These are the multiplicative counterparts of the parameter change
for the Hahn polynomials. The parameter range could well be enlarged.
\begin{align}
  &q^{\bm{\lambda}}=(a,b,q^N),\quad
  \bm{\delta}=(1,1,-1),\quad \kappa=q^{-1},\\
  &\epsilon=\Bigl\{\begin{array}{ll}
  1&\text{for }\ 0<a<1,\ \ 0<b<1,\\
  -1&\text{for }\ a>q^{1-N},\ \ b>q^{1-N},
  \end{array}\\
  &B(x\,;\bm{\lambda})=\epsilon(1-aq^x)(q^{x-N}-1),\quad
  D(x\,;\bm{\lambda})=\epsilon aq^{-1}(1-q^x)(q^{x-N}-b),\\
  &\mathcal{E}(n\,;\bm{\lambda})
  =\epsilon(q^{-n}-1)(1-abq^{n-1}),\quad
  \eta(x\,;\bm{\lambda})=q^{-x}-1,\\
  &P_n(\eta(x\,;\bm{\lambda})\,;\bm{\lambda})
  ={}_3\phi_2\Bigl(
  \genfrac{}{}{0pt}{}{q^{-n},\,abq^{n-1},\,q^{-x}}
  {a,\,q^{-N}}\Bigm|q\,;q\Bigr)
  =Q_n(1+\eta(x\,;\bm{\lambda})\,;aq^{-1},bq^{-1},N|q),\\
  &\phi_0(x\,;\bm{\lambda})^2
  =\frac{(q\,;q)_N}{(q\,;q)_x\,(q\,;q)_{N-x}}\,
  \frac{(a;q)_x\,(b\,;q)_{N-x}}{(b\,;q)_N\,a^x}\,,\\
  &d_n(\bm{\lambda})^2
  =\frac{(q\,;q)_N}{(q\,;q)_n\,(q\,;q)_{N-n}}\,
  \frac{(a,abq^{-1};q)_n}{(abq^N,b\,;q)_n\,a^n}\,
  \frac{1-abq^{2n-1}}{1-abq^{-1}}
  \times\frac{(b\,;q)_N\,a^N}{(ab\,;q)_N},\\
  &A_n(\bm{\lambda})
  =-\frac{(q^{n-N}-1)(1-aq^n)(1-abq^{n-1})}{(1-abq^{2n-1})(1-abq^{2n})}\,,\\
  &C_n(\bm{\lambda})
  =-aq^{n-N-1}\,
  \frac{(1-q^n)(1-abq^{n+N-1})(1-bq^{n-1})}{(1-abq^{2n-2})(1-abq^{2n-1})},\\
  &R_1(z\,;\bm{\lambda})=(q^{-\frac12}-q^{\frac12})^2z',\quad
  z'\eqdef z+\epsilon(1+abq^{-1}),\\
  &R_0(z\,;\bm{\lambda})=(q^{-\frac12}-q^{\frac12})^2
  \bigl(z^{\prime\,2}-ab(1+q^{-1})^2\bigr),\\
  &R_{-1}(z\,;\bm{\lambda})=(q^{-\frac12}-q^{\frac12})^2
  \Bigl(z^{\prime\,2}
  -\bigl(a(1+bq^{-1})+(1+aq^{-1})q^{-N}\bigr)\epsilon z'\n
  &\qquad\qquad\qquad\qquad\qquad\qquad
  +a(1+q^{-1})\bigl((a-1)bq^{-1}+(1+bq^{-1})q^{-N}\bigr)\Bigr),\\
  &B(x\,;\bm{\lambda}+\bm{\delta})=B(x+1\,;\bm{\lambda}),\quad
  D(x\,;\bm{\lambda}+\bm{\delta})=q^2D(x\,;\bm{\lambda}),\\
  &\varphi(x\,;\bm{\lambda})=q^{-x},\quad
  f_n(\bm{\lambda})=\mathcal{E}(n\,;\bm{\lambda}),\quad
  b_n(\bm{\lambda})=1.
\end{align}
Obviously the $q$-Hahn and dual $q$-Hahn are dual to each other.

\subsubsection{dual $q$-Hahn [KS3.7] (dual to \S\ref{[KS3.6]})}
\label{[KS3.7]}

For obvious reasons, we adopt the same parameters $(a,b)$ for
the $q$-Hahn and dual $q$-Hahn polynomials:
\begin{align}
  &q^{\bm{\lambda}}=(a,b,q^N),\quad
  \bm{\delta}=(1,0,-1),\quad \kappa=q^{-1},\\
  &\epsilon=\Bigl\{\begin{array}{ll}
  1&\text{for }\ 0<a<1,\ \ 0<b<1,\\
  -1&\text{for }\ a>q^{1-N},\ \ b>q^{1-N},
  \end{array}\\
  &B(x\,;\bm{\lambda})=
  \frac{(q^{x-N}-1)(1-aq^x)(1-abq^{x-1})}
  {(1-abq^{2x-1})(1-abq^{2x})},\\
  &D(x\,;\bm{\lambda})=aq^{x-N-1}
  \frac{(1-q^x)(1-abq^{x+N-1})(1-bq^{x-1})}
  {(1-abq^{2x-2})(1-abq^{2x-1})},\\
  &\mathcal{E}(n\,;\bm{\lambda})=q^{-n}-1,\quad
  \eta(x\,;\bm{\lambda})=\epsilon(q^{-x}-1)(1-abq^{x-1}),\\
  &P_n(\eta(x\,;\bm{\lambda})\,;\bm{\lambda})
  ={}_3\phi_2\Bigl(
  \genfrac{}{}{0pt}{}{q^{-n},\,abq^{x-1},\,q^{-x}}
  {a,\,q^{-N}}\Bigm|q\,;q\Bigr)\n
  &\phantom{P_n(\eta(x\,;\bm{\lambda})\,;\bm{\lambda})}
  =R_n(\epsilon(1+abq^{-1}+\eta(x\,;\bm{\lambda}))\,;aq^{-1},bq^{-1},N|q),\\
  &\phi_0(x\,;\bm{\lambda})^2
  =\frac{(q\,;q)_N}{(q\,;q)_x\,(q\,;q)_{N-x}}\,
  \frac{(a,abq^{-1}\,;q)_x}{(abq^N,b\,;q)_x\,a^x}\,
  \frac{1-abq^{2x-1}}{1-abq^{-1}}\,,\\
  &d_n(\bm{\lambda})^2
  =\frac{(q\,;q)_N}{(q\,;q)_n\,(q\,;q)_{N-n}}\,
  \frac{(a\,;q)_n(b\,;q)_{N-n}}{(b;q)_N\,a^n}
  \times\frac{(b\,;q)_N\,a^N}{(ab;q)_N}\,,\\
  &A_n(\bm{\lambda})
  =-\epsilon(1-aq^n)(q^{n-N}-1),\quad
  C_n(\bm{\lambda})
  =-\epsilon aq^{-1}(1-q^n)(q^{n-N}-b),\\
  &R_1(z\,;\bm{\lambda})=(q^{-\frac12}-q^{\frac12})^2 z',\quad
  z'\eqdef z+1,\\
  &R_0(z\,;\bm{\lambda})=(q^{-\frac12}-q^{\frac12})^2
  z^{\prime\,2},\\
  &R_{-1}(z\,;\bm{\lambda})=(q^{-\frac12}-q^{\frac12})^2
  \Bigl(\epsilon(1+abq^{-1})z^{\prime\,2}
  -\epsilon\bigl(a(1+bq^{-1})+(1+aq^{-1})q^{-N}\bigr)z'\n
  &\qquad\qquad\qquad\qquad\qquad\quad
  +\epsilon a(1+q^{-1})q^{-N}\Bigr),\\
  &B(x\,;\bm{\lambda}+\bm{\delta})=B(x+1\,;\bm{\lambda})\,
  \frac{1-abq^{2x+2}}{1-abq^{2x}}\,,\quad
  D(x\,;\bm{\lambda}+\bm{\delta})=q^2D(x\,;\bm{\lambda})\,
  \frac{1-abq^{2x-2}}{1-abq^{2x}}\,,\\
  &\varphi(x\,;\bm{\lambda})=\frac{q^{-x}-abq^x}{1-ab},\quad
  f_n(\bm{\lambda})=\mathcal{E}(n\,;\bm{\lambda}),\quad
  b_n(\bm{\lambda})=1.
\end{align}

\subsubsection{quantum $q$-Krawtchouk [KS3.14]}
\label{[KS3.14]}

\begin{align}
  &q^{\bm{\lambda}}=(p,q^N),\quad
  \bm{\delta}=(1,-1),\quad \kappa=q,\quad p>q^{-N},\\
  &B(x\,;\bm{\lambda})=p^{-1}q^x(q^{x-N}-1),\quad
  D(x\,;\bm{\lambda})=(1-q^x)(1-p^{-1}q^{x-N-1}),\\
  &\mathcal{E}(n\,;\bm{\lambda})=1-q^n,\quad
  \eta(x\,;\bm{\lambda})=q^{-x}-1,\\
  &P_n(\eta(x\,;\bm{\lambda})\,;\bm{\lambda})
  ={}_2\phi_1\Bigl(
  \genfrac{}{}{0pt}{}{q^{-n},\,q^{-x}}{q^{-N}}\Bigm|q\,;pq^{n+1}\Bigr)
  =K^{\text{qtm}}_n(1+\eta(x\,;\bm{\lambda})\,;p,N\,;q),\\
  &\phi_0(x\,;\bm{\lambda})^2
  =\frac{(q\,;q)_N}{(q\,;q)_x(q\,;q)_{N-x}}\,
  \frac{p^{-x}q^{x(x-1-N)}}{(p^{-1}q^{-N}\,;q)_x}\,,\\
  &d_n(\bm{\lambda})^2
  =\frac{(q\,;q)_N}{(q\,;q)_n(q\,;q)_{N-n}}\,
  \frac{p^{-n}q^{-Nn}}{(p^{-1}q^{-n}\,;q)_n}\,
  \times(p^{-1}q^{-N}\,;q)_N,\\
  &A_n(\bm{\lambda})=-p^{-1}q^{-n-N-1}(1-q^{N-n}),\quad
  C_n(\bm{\lambda})=-(q^{-n}-1)(1-p^{-1}q^{-n}),\\
  &R_1(z\,;\bm{\lambda})=(q^{-\frac12}-q^{\frac12})^2 z',\quad
  z'\eqdef z-1,\\
  &R_0(z\,;\bm{\lambda})=(q^{-\frac12}-q^{\frac12})^2
  z^{\prime\,2},\\
  &R_{-1}(z\,;\bm{\lambda})=(q^{-\frac12}-q^{\frac12})^2
  \bigl(z^{\prime\,2}+p^{-1}(1+p+q^{-N-1})z'
  +p^{-1}(1+q^{-1})\bigr),\\
  &B(x\,;\bm{\lambda}+\bm{\delta})=q^{-1}B(x+1\,;\bm{\lambda}),\quad
  D(x\,;\bm{\lambda}+\bm{\delta})=qD(x\,;\bm{\lambda}),\\
  &\varphi(x\,;\bm{\lambda})=q^{-x},\quad
  f_n(\bm{\lambda})=\mathcal{E}(n\,;\bm{\lambda}),\quad
  b_n(\bm{\lambda})=1.
\end{align}

The dual quantum $q$-Krawtchouk polynomial with the same parameter and
\begin{align}
  &B(x\,;\bm{\lambda})=p^{-1}q^{-x-N-1}(1-q^{N-x}),\quad
  D(x\,;\bm{\lambda})=(q^{-x}-1)(1-p^{-1}q^{-x}),
  \label{dualqqkrawt1}\\
  &\mathcal{E}(n\,;\bm{\lambda})=q^{-n}-1,\quad
  \eta(x\,;\bm{\lambda})=1-q^{x},\\
  &A_n(\bm{\lambda})=-p^{-1}q^n(q^{n-N}-1),\quad
  C_n(\bm{\lambda})=-(1-q^n)(1-p^{-1}q^{n-N-1}),
  \label{dualqqkrawt3}
\end{align}
has not been reported in Koekoek-Swarttouw \cite{koeswart}.
It is interesting to note that the functions $B(x)$ and $D(x)$ in
\eqref{dualqqkrawt1} are related to those of the affine
$q$-Krawtchouk polynomial \S\ref{[KS3.16]} with the change of
variable $x\leftrightarrow N-x$, $B\leftrightarrow D$:
\begin{equation}
  B(x)=D^{\text{affine $q$-Krawtchouk}}(N-x),\quad
  D(x)=B^{\text{affine $q$-Krawtchouk}}(N-x),
\end{equation}
with the change of the parameter
$p^{\text{affine $q$-Krawtchouk}}\to p^{-1}q^{-N-1}$.

\subsubsection{$q$-Krawtchouk [KS3.15] and dual
$q$-Krawtchouk [KS3.17]}
\label{[KS3.15]}

\begin{align}
  &q^{\bm{\lambda}}=(p,q^N),\
  \bm{\delta}=(2,-1),\quad \kappa=q^{-1},\quad p>0,\\
  &B(x\,;\bm{\lambda})=q^{x-N}-1,\quad
  D(x\,;\bm{\lambda})=p(1-q^x),\\
  &\mathcal{E}(n\,;\bm{\lambda})=(q^{-n}-1)(1+pq^n),\quad
  \eta(x\,;\bm{\lambda})=q^{-x}-1,\\
  &P_n(\eta(x\,;\bm{\lambda})\,;\bm{\lambda})
  ={}_3\phi_2\Bigl(
  \genfrac{}{}{0pt}{}{q^{-n},\,q^{-x},\,-pq^n}{q^{-N},\,0}\Bigm|q\,;q\Bigr)
  =K_n(1+\eta(x\,;\bm{\lambda})\,;p,N\,;q),\\
  &\phi_0(x\,;\bm{\lambda})^2=\frac{(q\,;q)_N}{(q\,;q)_x(q\,;q)_{N-x}}\,
  p^{-x}q^{\frac12x(x-1)-xN},\\
  &d_n(\bm{\lambda})^2
  =\frac{(q\,;q)_N}{(q;q)_n(q;q)_{N-n}}\,
  \frac{(-p\,;q)_n}{(-pq^{N+1}\,;q)_n\,p^nq^{\frac12n(n+1)}}\,
  \frac{1+pq^{2n}}{1+p}
  \times\frac{p^{N}q^{\frac12N(N+1)}}{(-pq\,;q)_N},\\
  &A_n(\bm{\lambda})=-\frac{(q^{n-N}-1)(1+pq^n)}{(1+pq^{2n})(1+pq^{2n+1})}\,,
  \quad
  C_n(\bm{\lambda})=-pq^{2n-N-1}\frac{(1-q^n)(1+pq^{n+N})}
  {(1+pq^{2n-1})(1+pq^{2n})}\,,\\
  &R_1(z\,;\bm{\lambda})=(q^{-\frac12}-q^{\frac12})^2 z',\quad
  z'\eqdef z+1-p,\\
  &R_0(z\,;\bm{\lambda})=(q^{-\frac12}-q^{\frac12})^2
  \bigl(z^{\prime\,2}+p(q^{-\frac12}+q^{\frac12})^2\bigr),\\
  &R_{-1}(z\,;\bm{\lambda})=(q^{-\frac12}-q^{\frac12})^2
  \bigl(z^{\prime\,2}+(p-q^{-N})z'
  +p(1+q^{-1})(1-q^{-N})\bigr),\\
  &B(x\,;\bm{\lambda}+\bm{\delta})=B(x+1\,;\bm{\lambda}),\quad
  D(x\,;\bm{\lambda}+\bm{\delta})=q^2D(x\,;\bm{\lambda}),\\
  &\varphi(x\,;\bm{\lambda})=q^{-x},\quad
  f_n(\bm{\lambda})=\mathcal{E}(n\,;\bm{\lambda}),\quad
  b_n(\bm{\lambda})=1.
\end{align}

The dual $q$-Krawtchouk polynomial in this parametrisation has
\begin{align}
  &B(x\,;\bm{\lambda})=\frac{(q^{x-N}-1)(1+pq^x)}{(1+pq^{2x})(1+pq^{2x+1})},
  \quad
  D(x\,;\bm{\lambda})=pq^{2x-N-1}
  \frac{(1-q^x)(1+pq^{x+N})}{(1+pq^{2x-1})(1+pq^{2x})},\\
  &\mathcal{E}(n\,;\bm{\lambda})=q^{-n}-1,\quad
  \eta(x\,;\bm{\lambda})=(q^{-x}-1)(1+pq^x),\\
  &A_n(\bm{\lambda})=-(q^{n-N}-1),\quad
  C_n(\bm{\lambda})=-p(1-q^n).
\end{align}
These are to be compared with the standard parametrisation of
the dual $q$-Krawtchouk [KS3.17], in which the parameter $c$ is
identified as
\begin{equation}
  c=-pq^N.
\end{equation}

\paragraph{dual $q$-Krawtchouk [KS3.17] in the standard parametrisation}

\begin{align}
  &q^{\bm{\lambda}}=(c,q^N),\quad
  \bm{\delta}=(0,-1),\quad \kappa=q^{-1},\quad c<0,\\
  &B(x\,;\bm{\lambda})=\frac{(q^{x-N}-1)(1-cq^{x-N})}
  {(1-cq^{2x-N})(1-cq^{2x+1-N})},\\
  &D(x\,;\bm{\lambda})=-cq^{2x-2N-1}\frac{(1-q^x)(1-cq^x)}
  {(1-cq^{2x-1-N})(1-cq^{2x-N})},\\
  &\mathcal{E}(n\,;\bm{\lambda})=q^{-n}-1,\quad
  \eta(x\,;\bm{\lambda})=(q^{-x}-1)(1-cq^{x-N}),\\
  &P_n(\eta(x\,;\bm{\lambda})\,;\bm{\lambda})
  ={}_3\phi_2\Bigl(
  \genfrac{}{}{0pt}{}{q^{-n},\,q^{-x},\,cq^{x-N}}{q^{-N},\,0}\Bigm|q\,;q\Bigr)
  =K_n(1+cq^{-N}+\eta(x\,;\bm{\lambda})\,;c,N|q),\\
  &\phi_0(x\,;\bm{\lambda})^2
  =\frac{(q\,;q)_N}{(q\,;q)_x(q\,;q)_{N-x}}\,
  \frac{(cq^{-N}\,;q)_x\,q^{Nx-\frac12x(x+1)}}{(cq\,;q)_x\,(-c)^x}\,
  \frac{1-cq^{2x-N}}{1-cq^{-N}}\,,\\
  &d_n(\bm{\lambda})^2
  =\frac{(q\,;q)_N}{(q\,;q)_n(q\,;q)_{N-n}}\,
  (-c)^{-n}q^{\frac12n(n-1)}
  \times\frac{1}{(c^{-1}\,;q)_N},\\
  &A_n(\bm{\lambda})=-(q^{n-N}-1),\quad
  C_n(\bm{\lambda})=cq^{-N}(1-q^n),\\
  &R_1(z\,;\bm{\lambda})=(q^{-\frac12}-q^{\frac12})^2 z',\quad
  z'\eqdef z+1,\\
  &R_0(z\,;\bm{\lambda})=(q^{-\frac12}-q^{\frac12})^2
  z^{\prime\,2},\\
  &R_{-1}(z\,;\bm{\lambda})=(q^{-\frac12}-q^{\frac12})^2
  \bigl((1+cq^{-N})z^{\prime\,2}-(1+c)q^{-N}z'\bigr),\\
  &B(x\,;\bm{\lambda}+\bm{\delta})=B(x+1\,;\bm{\lambda})\,
  \frac{1-cq^{2x+3-N}}{1-cq^{2x+1-N}},\,\,
  D(x\,;\bm{\lambda}+\bm{\delta})=q^2D(x\,;\bm{\lambda})\,
  \frac{1-cq^{2x-1-N}}{1-cq^{2x+1-N}},\\
  &\varphi(x\,;\bm{\lambda})=\frac{q^{-x}-cq^{1-N}q^x}{1-cq^{1-N}},\quad
  f_n(\bm{\lambda})=\mathcal{E}(n\,;\bm{\lambda}),\quad
  b_n(\bm{\lambda})=1.
\end{align}

\subsubsection{affine $q$-Krawtchouk [KS3.16] (self-dual)}
\label{[KS3.16]}

\begin{align}
  &q^{\bm{\lambda}}=(p,q^N),\quad
  \bm{\delta}=(1,-1),\quad \kappa=q^{-1},\quad 0<p<q^{-1}, \\
  &B(x\,;\bm{\lambda})=(q^{x-N}-1)(1-pq^{x+1}),\quad
  D(x\,;\bm{\lambda})=pq^{x-N}(1-q^x),\\
  &\mathcal{E}(n\,;\bm{\lambda})=q^{-n}-1,\quad
  \eta(x\,;\bm{\lambda})=q^{-x}-1,\\
  &P_n(\eta(x\,;\bm{\lambda})\,;\bm{\lambda})
  ={}_3\phi_2\Bigl(
  \genfrac{}{}{0pt}{}{q^{-n},\,q^{-x},\,0}{pq,\,q^{-N}}\Bigm|q\,;q\Bigr)
  =K^{\text{aff}}_n(1+\eta(x\,;\bm{\lambda})\,;p,N\,;q),\\
  &\phi_0(x\,;\bm{\lambda})^2=\frac{(q\,;q)_N}{(q\,;q)_x(q\,;q)_{N-x}}\,
  \frac{(pq\,;q)_x}{(pq)^x}\,,\\
  &d_n(\bm{\lambda})^2
  =\frac{(q\,;q)_N}{(q\,;q)_n(q\,;q)_{N-n}}\,
  \frac{(pq\,;q)_n}{(pq)^n}\times(pq)^N,\\
  &A_n(\bm{\lambda})=-(q^{n-N}-1)(1-pq^{n+1}),\quad
  C_n(\bm{\lambda})=-pq^{n-N}(1-q^n),\\
  &R_1(z\,;\bm{\lambda})=(q^{-\frac12}-q^{\frac12})^2 z',\quad
  z'\eqdef z+1,\\
  &R_0(z\,;\bm{\lambda})=(q^{-\frac12}-q^{\frac12})^2
  z^{\prime\,2},\\
  &R_{-1}(z\,;\bm{\lambda})=(q^{-\frac12}-q^{\frac12})^2
  \bigl(z^{\prime\,2}-(pq+(1+p)q^{-N})z'
  +p(1+q)q^{-N}\bigr),\\
  &B(x\,;\bm{\lambda}+\bm{\delta})=B(x+1\,;\bm{\lambda}),\quad
  D(x\,;\bm{\lambda}+\bm{\delta})=q^2D(x\,;\bm{\lambda}),\\
  &\varphi(x\,;\bm{\lambda})=q^{-x},\quad
  f_n(\bm{\lambda})=\mathcal{E}(n\,;\bm{\lambda}),\quad
  b_n(\bm{\lambda})=1.
\end{align}

\subsection{Infinite Dimensional Case (II)}

In this subsection $x$ takes an infinite range of values
\[
  x=0,1,2,\ldots,.
\]
In contrast to the finite dimensional case, the structure of the
polynomials is severely constrained by the asymptotic forms of
the functions $B(x)$ and $D(x)$, which are determined in Appendix A.
It is easy to see for $\eta(x)=x$ (\romannumeral1) \eqref{etaform1}, 
a quadratic spectrum
$\mathcal{E}(n)\sim n(n+\alpha)$ is not possible.
In order this to happen in \eqref{lowtri}, $B(x)$ and $D(x)$ must
have the same coefficient $\beta$ with the different sign for the
leading quadratic term, $B(x)\sim \beta x^2$, $D(x)\sim -\beta x^2$,
which violates the positivity of $B(x)$ and $D(x)$.
It is also easy to see $\eta(x)=\epsilon'x(x+d)$ (\romannumeral2) 
\eqref{etaform2} is not
possible. If the leading power of $B(x)$ and $D(x)$ is quadratic,
the coefficient must be the same with the same sign,
$B(x)\sim \beta x^2$, $D(x)\sim \beta x^2$, $\beta>0$.
Then for large $n$ \eqref{lowtri} gives a negative leading term
\begin{equation}
  \widetilde{\mathcal{H}}\eta(x)^n\sim -\beta n^2 \eta(x)^n+\cdots,
\end{equation}
which cannot give a normalisable eigenvector for a positive definite 
Hamiltonian.
If the leading power of $B(x)$ and $D(x)$ is linear, the results
in Appendix A \eqref{Btilderes} and \eqref{Dtilderes} tell that
$r^{(0)}_1=0$ and $B(x)$ and $D(x)$ have the opposite sign leading terms
\begin{align}
  B(x)&=\frac{r_{-1}^{(0)}\eta(x)}
  {\eta(1)B(0)\bigl(\eta(x+1)-\eta(x-1)\bigr)}
  +\text{lower order},\\
  D(x)&=-\frac{r_{-1}^{(0)}\eta(x)}
  {\eta(1)B(0)\bigl(\eta(x+1)-\eta(x-1)\bigr)}
  +\text{lower order}.
\end{align}
This simply contradicts the positivity of $B(x)$ and $D(x)$.
If the leading power of $B(x)$ and $D(x)$ is a constant or a negative
power in $x$, then eigenpolynomials in $\eta(x)$ of
$\widetilde{\mathcal H}$ do not exist.
Thus we have two self-dual polynomials in $x$, the Meixner \S\ref{[KS1.9]}
and Charlier \S\ref{[KS1.12]}.

As for the $q$-polynomials, the consequences of the asymptotic
behaviours of $B(x)$ and  $D(x)$ are easy to see.
If these functions are {\em bounded\/} as in the cases of the
sinusoidal coordinates $\eta(x)=q^{-x}-1$ (\romannumeral4) 
\eqref{etaform4}, and
$\eta(x)=\epsilon'(q^{-x}-1)(1-dq^x)$ (\romannumeral5) \eqref{etaform5},
the Hamiltonians are {\em bounded\/} and eigenvalues cannot take
{\em unbounded\/} forms $\mathcal{E}(n)=q^{-n}-1$,
$(q^{-n}-1)(1-\alpha q^n)$ and vice versa.
As mentioned in Appendix A, there is no self-dual $q$-polynomial of
infinite dimension.

\subsubsection{Meixner [KS1.9] (self-dual)}
\label{[KS1.9]}

This is the first example of the self-dual polynomial \eqref{MeixnerP},
which is symmetric under interchange $x\leftrightarrow n$, with
$\mathcal{E}(n)=n$ and $\eta(x)=x$ \eqref{MeixnerEeta}.
Compare \eqref{MeixnerBD}$\leftrightarrow$\eqref{MeixnerAC},
\eqref{Meixnerphi0d}:
\begin{align}
  &\bm{\lambda}=(\beta,c),\quad
  \bm{\delta}=(1,0),\quad \kappa=1,\quad \beta>0,\quad 0<c<1,\\
  &B(x\,;\bm{\lambda})=\frac{c}{1-c}(x+\beta),\quad
  D(x\,;\bm{\lambda})=\frac{1}{1-c}\,x,
  \label{MeixnerBD}\\
  &\mathcal{E}(n\,;\bm{\lambda})=n,\quad
  \eta(x\,;\bm{\lambda})=x,
  \label{MeixnerEeta}\\
  &P_n(\eta(x\,;\bm{\lambda})\,;\bm{\lambda})
  ={}_2F_1\Bigl(
  \genfrac{}{}{0pt}{}{-n,\,-x}{\beta}\Bigm|1-c^{-1}\Bigr)
  =M_n(\eta(x\,;\bm{\lambda})\,;\beta,c),
  \label{MeixnerP}\\
  &\phi_0(x\,;\bm{\lambda})^2=\frac{(\beta)_x\,c^x}{x!}\,,\quad
  d_n(\bm{\lambda})^2
  =\frac{(\beta)_n\,c^n}{n!}\times(1-c)^{\beta},
  \label{Meixnerphi0d}\\
  &A_n(\bm{\lambda})=-\frac{c}{1-c}(n+\beta),\quad
  C_n(\bm{\lambda})=-\frac{1}{1-c}\,n,
  \label{MeixnerAC}\\
  &R_1(z\,;\bm{\lambda})=0,\quad
  R_0(z\,;\bm{\lambda})=1,\quad
  R_{-1}(z\,;\bm{\lambda})=-\frac{1+c}{1-c}\,z-\frac{\beta c}{1-c},\\
  &B(x\,;\bm{\lambda}+\bm{\delta})=B(x+1\,;\bm{\lambda}),\quad
  D(x\,;\bm{\lambda}+\bm{\delta})=D(x\,;\bm{\lambda}),\\
  &\varphi(x\,;\bm{\lambda})=1,\quad
  f_n(\bm{\lambda})=\mathcal{E}(n\,;\bm{\lambda}),\quad
  b_n(\bm{\lambda})=1.
\end{align}

\subsubsection{Charlier [KS1.12] (self-dual)}
\label{[KS1.12]}

\begin{align}
  &\bm{\lambda}=a,\quad
  \bm{\delta}=0,\quad \kappa=1,\quad a>0,\\
  &B(x\,;\bm{\lambda})=a,\quad
  D(x\,;\bm{\lambda})=x,
  \label{charlBD}\\
  &\mathcal{E}(n\,;\bm{\lambda})=n,\quad
  \eta(x\,;\bm{\lambda})=x,
  \label{charlEeta}\\
  &P_n(\eta(x\,;\bm{\lambda})\,;\bm{\lambda})
  ={}_2F_0\Bigl(
  \genfrac{}{}{0pt}{}{-n,\,-x}{-}\Bigm|-a^{-1}\Bigr)
  =C_n(\eta(x\,;\bm{\lambda})\,;a),
  \label{charlP}\\
  &\phi_0(x\,;\bm{\lambda})^2=\frac{a^x}{x!}\,,\quad
  d_n(\bm{\lambda})^2
  =\frac{a^{n}}{n!}\times e^{-a},
  \label{charlphi0d}\\
  &A_n(\bm{\lambda})=-a,\quad
  C_n(\bm{\lambda})=-n,
  \label{charlAC}\\
  &R_1(z\,;\bm{\lambda})=0,\quad
  R_0(z\,;\bm{\lambda})=1,\quad
  R_{-1}(z\,;\bm{\lambda})=-z-a,\\
  &B(x\,;\bm{\lambda}+\bm{\delta})=B(x+1\,;\bm{\lambda}),\quad
  D(x\,;\bm{\lambda}+\bm{\delta})=D(x\,;\bm{\lambda}),\\
  &\varphi(x\,;\bm{\lambda})=1,\quad
  f_n(\bm{\lambda})=\mathcal{E}(n\,;\bm{\lambda}),\quad
  b_n(\bm{\lambda})=1.
\end{align}

\subsubsection{little $q$-Jacobi [KS3.12]}
\label{[KS3.12]}

The universal normalisation \eqref{Pzero} differs from the standard
normalisation of the little $q$-Jacobi polynomial $p_n$ as shown
explicitly in \eqref{littleqjacobinorm}. This does not affect the
normalisation measure $\phi_0(x)^2$ \eqref{littleqjacobiphi0} but
the normalisation constants $d_n^2$ \eqref{littleqjacobidn} and the
coefficients of the three term recurrence $A_n$ and $C_n$
\eqref{littleqjacobiAn}--\eqref{littleqjacobiCn} are different from
the standard ones:
\begin{align}
  &q^{\bm{\lambda}}=(a,b),\quad
  \bm{\delta}=(1,1),\quad \kappa=q^{-1},\quad 0<a<q^{-1},\quad b<q^{-1},\\
  &B(x\,;\bm{\lambda})=a(q^{-x}-bq),\quad
  D(x\,;\bm{\lambda})=q^{-x}-1,\\
  &\mathcal{E}(n\,;\bm{\lambda})=(q^{-n}-1)(1-abq^{n+1}),\quad
  \eta(x\,;\bm{\lambda})=1-q^x,\\
  &P_n(\eta(x\,;\bm{\lambda})\,;\bm{\lambda})
  =(-a)^{-n}q^{-\frac12n(n+1)}\frac{(aq\,;q)_n}{(bq\,;q)_n}\,
  {}_2\phi_1\Bigl(
  \genfrac{}{}{0pt}{}{q^{-n},\,abq^{n+1}}{aq}\Bigm|q\,;q^{x+1}\Bigr)\n
  &\phantom{P_n(\eta(x\,;\bm{\lambda})\,;\bm{\lambda})}
  =(-a)^{-n}q^{-\frac12n(n+1)}\frac{(aq\,;q)_n}{(bq\,;q)_n}\,
  p_n(1-\eta(x\,;\bm{\lambda})\,;a,b|q),
  \label{littleqjacobinorm}\\
  &\phi_0(x\,;\bm{\lambda})^2=\frac{(bq\,;q)_x}{(q\,;q)_x}(aq)^x,
  \label{littleqjacobiphi0}\\
  &d_n(\bm{\lambda})^2
  =\frac{(bq,abq\,;q)_n\,a^nq^{n^2}}{(q,aq\,;q)_n}\,
  \frac{1-abq^{2n+1}}{1-abq}
  \times\frac{(aq\,;q)_{\infty}}{(abq^2\,;q)_{\infty}}\,,
  \label{littleqjacobidn}\\
  &A_n(\bm{\lambda})=-aq^{2n+1}\frac{(1-bq^{n+1})(1-abq^{n+1})}
  {(1-abq^{2n+1})(1-abq^{2n+2})}\,,
  \label{littleqjacobiAn}\\
  &C_n(\bm{\lambda})=-\frac{(1-q^n)(1-aq^n)}
  {(1-abq^{2n})(1-abq^{2n+1})}\,,
  \label{littleqjacobiCn}\\
  &R_1(z\,;\bm{\lambda})=(q^{-\frac12}-q^{\frac12})^2 z',\quad
  z'\eqdef z+1+abq,\\
  &R_0(z\,;\bm{\lambda})=(q^{-\frac12}-q^{\frac12})^2
  \bigl(z^{\prime\,2}-ab(1+q)^2\bigr),\\
  &R_{-1}(z\,;\bm{\lambda})=(q^{-\frac12}-q^{\frac12})^2
  \bigl(-z^{\prime\,2}+(1+a)z'-a(1+q)(1-bq)\bigr),\\
  &B(x\,;\bm{\lambda}+\bm{\delta})=q^2B(x+1\,;\bm{\lambda}),\quad
  D(x\,;\bm{\lambda}+\bm{\delta})=D(x\,;\bm{\lambda}),\\
  &\varphi(x\,;\bm{\lambda})=q^x,\quad
  f_n(\bm{\lambda})=\mathcal{E}(n\,;\bm{\lambda}),\quad
  b_n(\bm{\lambda})=1.
\end{align}

The dual little $q$-Jacobi polynomial with the same parameter and
\begin{align}
  &B(x\,;\bm{\lambda})=aq^{2x+1}\frac{(1-bq^{x+1})(1-abq^{x+1})}
  {(1-abq^{2x+1})(1-abq^{2x+2})},\ \,
  D(x\,;\bm{\lambda})=\frac{(1-q^x)(1-aq^x)}
  {(1-abq^{2x})(1-abq^{2x+1})},\\
  &\mathcal{E}(n\,;\bm{\lambda})=1-q^n,\quad
  \eta(x\,;\bm{\lambda})=(q^{-x}-1)(1-abq^{x+1}),\\
  &A_n(\bm{\lambda})=-a(q^{-n}-bq),\quad
  C_n(\bm{\lambda})=-(q^{-n}-1),
\end{align}
has not been reported in Koekoek-Swarttouw \cite{koeswart}.
The dual little $q$-Jacobi polynomial was introduced by
Atakishiyev and Klimyk \cite{atakishi2}, which has the same three
term recurrence relation and the orthogonality measure as above.
The difference equation was not mentioned there.

\subsubsection{$q$-Meixner [KS3.13]}
\label{[KS3.13]}

\begin{align}
  &q^{\bm{\lambda}}=(b,c),\quad
  \bm{\delta}=(1,-1),\quad \kappa=q,\quad 0<b<q^{-1},\quad c>0,\\
  &B(x\,;\bm{\lambda})=cq^x(1-bq^{x+1}),\quad
  D(x\,;\bm{\lambda})=(1-q^x)(1+bcq^x),\\
  &\mathcal{E}(n\,;\bm{\lambda})=1-q^n,\quad
  \eta(x\,;\bm{\lambda})=q^{-x}-1,\\
  &P_n(\eta(x\,;\bm{\lambda})\,;\bm{\lambda})
  ={}_2\phi_1\Bigl(
  \genfrac{}{}{0pt}{}{q^{-n},\,q^{-x}}{bq}\Bigm|q\,;-c^{-1}q^{n+1}\Bigr)
  =M_n(1+\eta(x\,;\bm{\lambda})\,;b,c\,;q),\\
  &\phi_0(x\,;\bm{\lambda})^2=
  \frac{(bq\,;q)_x}{(q,-bcq\,;q)_x}\,c^xq^{\frac12x(x-1)},\quad
  d_n(\bm{\lambda})^2
  =\frac{(bq\,;q)_n}{(q,-c^{-1}q\,;q)_n}
  \times\frac{(-bcq\,;q)_{\infty}}{(-c\,;q)_{\infty}}\,,\\
  &A_n(\bm{\lambda})=-cq^{-n-1}(q^{-n}-bq),\quad
  C_n(\bm{\lambda})=-(q^{-n}-1)(1+cq^{-n}),\\
  &R_1(z\,;\bm{\lambda})=(q^{-\frac12}-q^{\frac12})^2 z',\quad
  z'\eqdef z-1,\\
  &R_0(z\,;\bm{\lambda})=(q^{-\frac12}-q^{\frac12})^2
  z^{\prime\,2},\\
  &R_{-1}(z\,;\bm{\lambda})=(q^{-\frac12}-q^{\frac12})^2
  \bigl(z^{\prime\,2}+(1-c-bc)z'-c(1+q^{-1})\bigr),\\
  &B(x\,;\bm{\lambda}+\bm{\delta})=q^{-2}B(x+1\,;\bm{\lambda}),\quad
  D(x\,;\bm{\lambda}+\bm{\delta})=D(x\,;\bm{\lambda}),\\
  &\varphi(x\,;\bm{\lambda})=q^{-x},\quad
  f_n(\bm{\lambda})=\mathcal{E}(n\,;\bm{\lambda}),\quad
  b_n(\bm{\lambda})=1.
\end{align}

The dual $q$-Meixner polynomial with the same parameter and
\begin{align}
  &B(x\,;\bm{\lambda})=cq^{-x-1}(q^{-x}-bq),\quad
  D(x\,;\bm{\lambda})=(q^{-x}-1)(1+cq^{-x}),
  \label{dualqmeixner1}\\
  &\mathcal{E}(n\,;\bm{\lambda})=q^{-n}-1,\quad
  \eta(x\,;\bm{\lambda})=1-q^x,\\
  &A_n(\bm{\lambda})=-cq^n(1-bq^{n+1}),\quad
  C_n(\bm{\lambda})=-(1-q^n)(1+bcq^n),
  \label{dualqmeixner3}
\end{align}
has not been reported in Koekoek-Swarttouw \cite{koeswart}.
This can be considered another $q$-version of the Meixner polynomial
\S\ref{[KS1.9]} with the sinusoidal coordinate $\eta(x)=1-q^x$.

\subsubsection{little $q$-Laguerre/Wall [KS3.20] (dual to \S\ref{[KS3.25]})}
\label{[KS3.20]}

The universal normalisation \eqref{Pzero} differs from the standard
normalisation of the little $q$-Laguerre/Wall polynomial $p_n$
as shown explicitly in \eqref{littleqlaguerrenorm}.
As other examples, $d_n^2$, $A_n$ and $C_n$ are different from the
standard ones:
\begin{align}
  &q^{\bm{\lambda}}=a,\quad
  \bm{\delta}=1,\quad \kappa=q^{-1},\quad 0<a<q^{-1},\\
  &B(x\,;\bm{\lambda})=aq^{-x},\quad
  D(x\,;\bm{\lambda})=q^{-x}-1,\\
  &\mathcal{E}(n\,;\bm{\lambda})=q^{-n}-1,\quad
  \eta(x\,;\bm{\lambda})=1-q^x,\\
  &P_n(\eta(x\,;\bm{\lambda})\,;\bm{\lambda})
  ={}_2\phi_0\Bigl(
  \genfrac{}{}{0pt}{}{q^{-n},\,q^{-x}}{-}\Bigm|q\,;a^{-1}q^x\Bigr)
  =(a^{-1}q^{-n}\,;q)_n\,p_n(1-\eta(x\,;\bm{\lambda})\,;a|q),
  \label{littleqlaguerrenorm}\\
  &\phi_0(x\,;\bm{\lambda})^2=\frac{(aq)^x}{(q\,;q)_x}\,,\quad
  d_n(\bm{\lambda})^2
  =\frac{a^nq^{n^2}}{(q,aq\,;q)_n}\times(aq\,;q)_{\infty}\,,\\
  &A_n(\bm{\lambda})=-aq^{2n+1},\quad
  C_n(\bm{\lambda})=-(1-q^n)(1-aq^n),\\
  &R_1(z\,;\bm{\lambda})=(q^{-\frac12}-q^{\frac12})^2 z',\quad
  z'\eqdef z+1,\\
  &R_0(z\,;\bm{\lambda})=(q^{-\frac12}-q^{\frac12})^2
  z^{\prime\,2},\\
  &R_{-1}(z\,;\bm{\lambda})=(q^{-\frac12}-q^{\frac12})^2
  \bigl(-z^{\prime\,2}+(1+a)z'-a(1+q)\bigr),\\
  &B(x\,;\bm{\lambda}+\bm{\delta})=q^2B(x+1\,;\bm{\lambda}),\quad
  D(x\,;\bm{\lambda}+\bm{\delta})=D(x\,;\bm{\lambda}),\\
  &\varphi(x\,;\bm{\lambda})=q^x,\quad
  f_n(\bm{\lambda})=\mathcal{E}(n\,;\bm{\lambda}),\quad
  b_n(\bm{\lambda})=1.
\end{align}

The dual little $q$-Laguerre/Wall polynomial with the same parameter
and
\begin{align}
  &B(x;\bm{\lambda})=aq^{2x+1},\quad
  D(x;\bm{\lambda})=(1-q^x)(1-aq^x),\\
  &\mathcal{E}(n\,;\bm{\lambda})=1-q^n,\quad
  \eta(x\,;\bm{\lambda})=q^{-x}-1,\\
  &A_n(\bm{\lambda})=-aq^{-n},\quad
  C_n(\bm{\lambda})=-(q^{-n}-1),
\end{align}
is Al-Salam Carlitz II polynomial as seen below.
This was reported by Atakishiyev and Klimyk \cite{atakishi2}.

\subsubsection{Al-Salam-Carlitz II [KS3.25] (dual to \S\ref{[KS3.20]})}
\label{[KS3.25]}

The universal normalisation \eqref{Pzero} differs from the standard
normalisation of the Al-Salam-Carlitz II polynomial $V_n^{(a)}$
as shown explicitly in \eqref{alsalamIInorm}.
As other examples, $d_n^2$, $A_n$ and $C_n$ are different from the
standard ones:
\begin{align}
  &q^{\bm{\lambda}}=a,\quad
  \bm{\delta}=0,\quad \kappa=q,\quad 0<a<q^{-1},\\
  &B(x\,;\bm{\lambda})=aq^{2x+1},\quad
  D(x\,;\bm{\lambda})=(1-q^x)(1-aq^x),\\
  &\mathcal{E}(n\,;\bm{\lambda})=1-q^n,\quad
  \eta(x\,;\bm{\lambda})=q^{-x}-1,\\
  &P_n(\eta(x\,;\bm{\lambda})\,;\bm{\lambda})
  ={}_2\phi_0\Bigl(
  \genfrac{}{}{0pt}{}{q^{-n},\,q^{-x}}{-}\Bigm|q\,;a^{-1}q^n\Bigr)
  =(-a^{-1}q^{\frac12(n-1)})^n\,V^{(a)}_n(1+\eta(x\,;\bm{\lambda})\,;q),
  \label{alsalamIInorm}\\
  &\phi_0(x\,;\bm{\lambda})^2=\frac{a^xq^{x^2}}{(q,aq\,;q)_x}\,,\quad
  d_n(\bm{\lambda})^2
  =\frac{(aq)^n}{(q\,;q)_n}\times(aq\,;q)_{\infty}\,,\\
  &A_n(\bm{\lambda})=-aq^{-n},\quad
  C_n(\bm{\lambda})=-(q^{-n}-1),\\
  &R_1(z\,;\bm{\lambda})=(q^{-\frac12}-q^{\frac12})^2 z',\quad
  z'\eqdef z-1,\\
  &R_0(z\,;\bm{\lambda})=(q^{-\frac12}-q^{\frac12})^2
  z^{\prime\,2},\\
  &R_{-1}(z\,;\bm{\lambda})=(q^{-\frac12}-q^{\frac12})^2
  \bigl(z^{\prime\,2}+(1+a)z'\bigr),\\
  &B(x\,;\bm{\lambda}+\bm{\delta})=q^{-2}B(x+1\,;\bm{\lambda}),\quad
  D(x\,;\bm{\lambda}+\bm{\delta})=D(x\,;\bm{\lambda}),\\
  &\varphi(x\,;\bm{\lambda})=q^{-x},\quad
  f_n(\bm{\lambda})=\mathcal{E}(n\,;\bm{\lambda}),\quad
  b_n(\bm{\lambda})=1.
\end{align}

\subsubsection{alternative $q$-Charlier [KS3.22]}
\label{[KS3.22]}

The universal normalisation \eqref{Pzero} differs from the standard
normalisation of the alternative $q$-Charlier polynomial $K_n$
as shown explicitly in \eqref{alcharliernorm}.
As other examples, $d_n^2$, $A_n$ and $C_n$ are different from the
standard ones:
\begin{align}
  &q^{\bm{\lambda}}=a,\quad
  \bm{\delta}=2,\quad \kappa=q^{-1},\quad a>0,\\
  &B(x\,;\bm{\lambda})=a,\quad
  D(x\,;\bm{\lambda})=q^{-x}-1,\\
  &\mathcal{E}(n\,;\bm{\lambda})=(q^{-n}-1)(1+aq^n),\quad
  \eta(x\,;\bm{\lambda})=1-q^x,\\
  &P_n(\eta(x\,;\bm{\lambda})\,;\bm{\lambda})
  =q^{nx}\,{}_2\phi_1\Bigl(
  \genfrac{}{}{0pt}{}{q^{-n},\,q^{-x}}{0}\Bigm|q\,;-a^{-1}q^{-n+1}\Bigr)\n
  &\phantom{P_n(\eta(x\,;\bm{\lambda})\,;\bm{\lambda})}
  =(-aq^n)^{-n}\,K_n(1-\eta(x\,;\bm{\lambda})\,;a\,;q),
  \label{alcharliernorm}\\
  &\phi_0(x\,;\bm{\lambda})^2=\frac{a^xq^{\frac12x(x+1)}}{(q\,;q)_x}\,,
  \ \,
  d_n(\bm{\lambda})^2
  =\frac{a^nq^{\frac12n(3n-1)}}{(q\,;q)_n}\,
  \frac{(-a\,;q)_{\infty}}{(-aq^n\,;q)_{\infty}}\,
  \frac{1+aq^{2n}}{1+a}
  \times\frac{1}{(-aq\,;q)_{\infty}}\,,\\
  &A_n(\bm{\lambda})=-aq^{3n+1}\frac{1+aq^n}{(1+aq^{2n})(1+aq^{2n+1})}\,,
  \quad
  C_n(\bm{\lambda})=-\frac{1-q^n}{(1+aq^{2n-1})(1+aq^{2n})}\,,\\
  &R_1(z\,;\bm{\lambda})=(q^{-\frac12}-q^{\frac12})^2 z',\quad
  z'\eqdef z+1-a,\\
  &R_0(z\,;\bm{\lambda})=(q^{-\frac12}-q^{\frac12})^2
  \bigl(z^{\prime\,2}+a(q^{-\frac12}+q^{\frac12})^2\bigr),\\
  &R_{-1}(z\,;\bm{\lambda})=(q^{-\frac12}-q^{\frac12})^2
  \bigl(-z^{\prime\,2}+z'-a(1+q)\bigr),\\
  &B(x\,;\bm{\lambda}+\bm{\delta})=q^2B(x+1\,;\bm{\lambda}),\quad
  D(x\,;\bm{\lambda}+\bm{\delta})=D(x\,;\bm{\lambda}),\\
  &\varphi(x\,;\bm{\lambda})=q^x,\quad
  f_n(\bm{\lambda})=\mathcal{E}(n\,;\bm{\lambda}),\quad
  b_n(\bm{\lambda})=1.
\end{align}

The dual alternative $q$-Charlier polynomial with the same parameter
and
\begin{align}
  &B(x\,;\bm{\lambda})=aq^{3x+1}\frac{1+aq^x}{(1+aq^{2x})(1+aq^{2x+1})},
  \quad
  D(x\,;\bm{\lambda})=\frac{1-q^x}{(1+aq^{2x-1})(1+aq^{2x})},\\
  &\mathcal{E}(n\,;\bm{\lambda})=1-q^n,\quad
  \eta(x\,;\bm{\lambda})=(q^{-x}-1)(1+aq^x),\\
  &A_n(\bm{\lambda})=-a,\quad
  C_n(\bm{\lambda})=-(q^{-n}-1),
\end{align}
has not been reported in Koekoek-Swarttouw \cite{koeswart}.
The dual alternative $q$-Charlier polynomial was introduced by
Atakishiyev and Klimyk \cite{atakishi2}, which has the same three
term recurrence relation and the orthogonality measure as above.
The difference equation was not mentioned there.

\subsubsection{$q$-Charlier [KS3.23]}
\label{[KS3.23]}

\begin{align}
  &q^{\bm{\lambda}}=a,\quad
  \bm{\delta}=-1,\quad \kappa=q,\quad a>0,\\
  &B(x\,;\bm{\lambda})=aq^x,\quad
  D(x\,;\bm{\lambda})=1-q^x,\\
  &\mathcal{E}(n\,;\bm{\lambda})=1-q^n,\quad
  \eta(x\,;\bm{\lambda})=q^{-x}-1,\\
  &P_n(\eta(x\,;\bm{\lambda})\,;\bm{\lambda})
  ={}_2\phi_1\Bigl(
  \genfrac{}{}{0pt}{}{q^{-n},\,q^{-x}}{0}\Bigm|q\,;-a^{-1}q^{n+1}\Bigr)
  =C_n(1+\eta(x\,;\bm{\lambda})\,;a\,;q),\\
  &\phi_0(x\,;\bm{\lambda})^2=\frac{a^xq^{\frac12x(x-1)}}{(q\,;q)_x}\,,\quad
  d_n(\bm{\lambda})^2
  =\frac{q^n}{(-a^{-1}q,q\,;q)_n}\times\frac{1}{(-a\,;q)_{\infty}}\,,\\
  &A_n(\bm{\lambda})=-aq^{-2n-1},\quad
  C_n(\bm{\lambda})=-(q^{-n}-1)(1+aq^{-n}),\\
  &R_1(z\,;\bm{\lambda})=(q^{-\frac12}-q^{\frac12})^2 z',\quad
  z'\eqdef z-1,\\
  &R_0(z\,;\bm{\lambda})=(q^{-\frac12}-q^{\frac12})^2
  z^{\prime\,2},\\
  &R_{-1}(z\,;\bm{\lambda})=(q^{-\frac12}-q^{\frac12})^2
  \bigl(z^{\prime\,2}+(1-a)z'-a(1+q^{-1})\bigr),\\
  &B(x\,;\bm{\lambda}+\bm{\delta})=q^{-2}B(x+1\,;\bm{\lambda}),\quad
  D(x\,;\bm{\lambda}+\bm{\delta})=D(x\,;\bm{\lambda}),\\
  &\varphi(x\,;\bm{\lambda})=q^{-x},\quad
  f_n(\bm{\lambda})=\mathcal{E}(n\,;\bm{\lambda}),\quad
  b_n(\bm{\lambda})=1.
\end{align}

The dual $q$-Charlier polynomial with the same parameter and
\begin{align}
  &B(x\,;\bm{\lambda})=aq^{-2x-1},\quad
  D(x\,;\bm{\lambda})=(q^{-x}-1)(1+aq^{-x}),
  \label{dualqcharlier1}\\
  &\mathcal{E}(n\,;\bm{\lambda})=q^{-n}-1,\quad
  \eta(x\,;\bm{\lambda})=1-q^x,\\
  &A_n(\bm{\lambda})=-aq^n,\quad
  C_n(\bm{\lambda})=-(1-q^n),
  \label{dualqcharlier3}
\end{align}
has not been reported in Koekoek-Swarttouw \cite{koeswart}.
This can be considered another $q$-version of the Charlier polynomial
\S\ref{[KS1.12]} with sinusoidal coordinate $\eta(x)=1-q^x$.
It should be emphasised that these two different $q$-versions of the
Charlier polynomial have very different characters.
The $q$-Charlier has a bounded spectrum and its dual has an unbounded
spectrum.
This is always the case for the infinite dimensional $q$-polynomials
and their duals.

\subsection{Other Polynomials}
\label{alternatives}

Here we will discuss orthogonal polynomials which are not listed
in Koekoek-Swarttouw's review \cite{koeswart}.
The naming of the polynomials is very tentative.

\subsubsection{alternative $q$-Hahn}
\label{aqHahn}

The $q$-Hahn polynomial \S\ref{[KS3.6]} is a $q$-version of the
Hahn polynomial \S\ref{[KS1.5]} with the sinusoidal coordinate
$\eta(x)=q^{-x}-1$.
The other $q$-version with $\eta(x)=1-q^{x}$ and the same spectrum
$\mathcal{E}(n)=\epsilon(q^{-n}-1)(1-abq^{n-1})$ as the $q$-Hahn
polynomial is as follows (the parameter range could well be enlarged):
\begin{align}
  &q^{\bm{\lambda}}=(a,b,q^N),\quad
  \bm{\delta}=(1,1,-1),\quad \kappa=q^{-1},\\
  &\epsilon=\Bigl\{\begin{array}{ll}
  1&\text{for }\ 0<a<1,\ \ 0<b<1,\\
  -1&\text{for }\ a>q^{1-N},\ \ b>q^{1-N},
  \end{array}\\
  &B(x\,;\bm{\lambda})=\epsilon aq^{-1}(1-q^{N-x})(q^{-x}-b),\quad
  D(x\,;\bm{\lambda})=\epsilon(1-aq^{N-x})(q^{-x}-1),\\
  &\mathcal{E}(n\,;\bm{\lambda})
  =\epsilon(q^{-n}-1)(1-abq^{n-1}),\quad
  \eta(x\,;\bm{\lambda})=1-q^x,\\
  &P_n(\eta(x\,;\bm{\lambda})\,;\bm{\lambda})
  ={}_3\phi_2\Bigl(
  \genfrac{}{}{0pt}{}{q^{-n},\,abq^{n-1},\,q^{-x}}
  {b,\,q^{-N}}\Bigm|q\,;a^{-1}q^{x+1-N}\Bigr),\\
  &\phi_0(x\,;\bm{\lambda})^2
  =\frac{(q\,;q)_N}{(q\,;q)_x\,(q\,;q)_{N-x}}\,
  \frac{a^x(a;q)_{N-x}\,(b\,;q)_x}{(a\,;q)_N}\,,\\
  &d_n(\bm{\lambda})^2
  =\frac{(q\,;q)_N}{(q\,;q)_n\,(q\,;q)_{N-n}}\,
  \frac{(b,abq^{-1};q)_n\,a^nq^{n(n-1)}}{(a,abq^N\,;q)_n}\,
  \frac{1-abq^{2n-1}}{1-abq^{-1}}
  \times\frac{(a\,;q)_N}{(ab\,;q)_N}\,,\\
  &A_n(\bm{\lambda})
  =-aq^{n+N}\,
  \frac{(q^{n-N}-1)(1-bq^n)(1-abq^{n-1})}{(1-abq^{2n-1})(1-abq^{2n})}\,,\\
  &C_n(\bm{\lambda})
  =-\frac{(1-q^n)(1-aq^{n-1})(1-abq^{n+N-1})}
  {(1-abq^{2n-2})(1-abq^{2n-1})}\,,\\
  &R_1(z\,;\bm{\lambda})=(q^{-\frac12}-q^{\frac12})^2 z',\quad
  z'\eqdef z+\epsilon(1+abq^{-1}),\\
  &R_0(z\,;\bm{\lambda})=(q^{-\frac12}-q^{\frac12})^2
  \bigl(z^{\prime\,2}-ab(1+q^{-1})^2\bigr),\\
  &R_{-1}(z\,;\bm{\lambda})=(q^{-\frac12}-q^{\frac12})^2
  \Bigl(-z^{\prime\,2}
  +\bigl(1+aq^{-1}+a(1+bq^{-1})q^N\bigr)\epsilon z'\n
  &\qquad\qquad\qquad\qquad\qquad\qquad
  -a(1+q^{-1})\bigl(1-b+b(1+aq^{-1})q^N\bigr)\Bigr),\\
  &B(x\,;\bm{\lambda}+\bm{\delta})=q^2B(x+1\,;\bm{\lambda}),\quad
  D(x\,;\bm{\lambda}+\bm{\delta})=D(x\,;\bm{\lambda}),\\
  &\varphi(x\,;\bm{\lambda})=q^x,\quad
  f_n(\bm{\lambda})=\mathcal{E}(n\,;\bm{\lambda}),\quad
  b_n(\bm{\lambda})=1.
\end{align}
This is the most generic form of the functions $B(x)$ and $D(x)$
for the sinusoidal coordinate $\eta(x)=1-q^x$ (quadratic polynomials
in $q^{-x}$) as mentioned \eqref{BDq3} in Appendix A.
It is interesting to note that the functions $B(x)$ and $D(x)$ are
related to those of the $q$-Hahn polynomial with the change of
variable $x\leftrightarrow N-x$, $B\leftrightarrow D$:
\begin{equation}
  B(x)=D^{\text{$q$-Hahn}}(N-x),\quad
  D(x)=B^{\text{$q$-Hahn}}(N-x).
\end{equation}

\subsubsection{alternative $q$-Krawtchouk}
\label{aqKrawt}

The $q$-Krawtchouk polynomial \S\ref{[KS3.15]} is a $q$-version of the
Krawtchouk polynomial \S\ref{[KS1.10]} with the sinusoidal coordinate
$\eta(x)=q^{-x}-1$.
The other $q$-version with $\eta(x)=1-q^{x}$ and the same spectrum
$\mathcal{E}(n)=(q^{-n}-1)(1+pq^n)$ is as follows:
\begin{align}
  &q^{\bm{\lambda}}=(p,q^N),\quad
  \bm{\delta}=(2,-1),\quad \kappa=q^{-1},\quad p>0,\\
  &B(x\,;\bm{\lambda})=p(1-q^{N-x}),\quad
  D(x\,;\bm{\lambda})=q^{-x}-1,\\
  &\mathcal{E}(n\,;\bm{\lambda})
  =(q^{-n}-1)(1+pq^n),\quad
  \eta(x\,;\bm{\lambda})=1-q^x,\\
  &P_n(\eta(x\,;\bm{\lambda})\,;\bm{\lambda})
  ={}_3\phi_1\Bigl(
  \genfrac{}{}{0pt}{}{q^{-n},\,-pq^n,\,q^{-x}}
  {q^{-N}}\Bigm|q\,;-p^{-1}q^{x-N}\Bigr),\\
  &\phi_0(x\,;\bm{\lambda})^2
  =\frac{(q\,;q)_N}{(q\,;q)_x\,(q\,;q)_{N-x}}\,
  p^xq^{\frac12x(x+1)},\\
  &d_n(\bm{\lambda})^2
  =\frac{(q\,;q)_N}{(q\,;q)_n\,(q\,;q)_{N-n}}\,
  \frac{(-p\,;q)_n\,p^nq^{\frac12n(3n-1)}}{(-pq^{N+1}\,;q)_n}\,
  \frac{1+pq^{2n}}{1+p}
  \times\frac{1}{(-pq\,;q)_N}\,,\\
  &A_n(\bm{\lambda})
  =-pq^{2n+N+1}\,
  \frac{(q^{n-N}-1)(1+pq^n)}{(1+pq^{2n})(1+pq^{2n+1})}\,,\quad
  C_n(\bm{\lambda})
  =-\frac{(1-q^n)(1+pq^{n+N})}{(1+pq^{2n-1})(1+pq^{2n})}\,,\\
  &R_1(z\,;\bm{\lambda})=(q^{-\frac12}-q^{\frac12})^2z',\quad
  z'\eqdef z+1-p,\\
  &R_0(z\,;\bm{\lambda})=(q^{-\frac12}-q^{\frac12})^2
  \bigl(z^{\prime\,2}+p(q^{-\frac12}+q^{\frac12})^2\bigr),\\
  &R_{-1}(z\,;\bm{\lambda})=(q^{-\frac12}-q^{\frac12})^2
  \bigl(-z^{\prime\,2}
  +(1-pq^N)z'-p(1+q)(1-q^N)\bigr),\\
  &B(x\,;\bm{\lambda}+\bm{\delta})=q^2B(x+1\,;\bm{\lambda}),\quad
  D(x\,;\bm{\lambda}+\bm{\delta})=D(x\,;\bm{\lambda}),\\
  &\varphi(x\,;\bm{\lambda})=q^x,\quad
  f_n(\bm{\lambda})=\mathcal{E}(n\,;\bm{\lambda}),\quad
  b_n(\bm{\lambda})=1.
\end{align}
The functions $B(x)$ and $D(x)$ are related to those of the
$q$-Krawtchouk polynomial with the change of variable
$x\leftrightarrow N-x$, $B\leftrightarrow D$:
\begin{equation}
  B(x)=D^{\text{$q$-Krawtchouk}}(N-x),\quad
  D(x)=B^{\text{$q$-Krawtchouk}}(N-x).
\end{equation}

\subsubsection{alternative affine $q$-Krawtchouk (self-dual)}
\label{aafqKrawt}

The affine $q$-Krawtchouk \S\ref{[KS3.16]} is a self-dual
polynomial with $\eta(x)=q^{-x}-1$ and $\mathcal{E}(n)=q^{-n}-1$.
There is another self-dual polynomial with $\eta(x)=1-q^{x}$ and
$\mathcal{E}(n)=1-q^{n}$:
\begin{align}
  &q^{\bm{\lambda}}=(p,q^N),\quad
  \bm{\delta}=(1,-1),\quad \kappa=q,\quad p>q^{-N}, \\
  &B(x\,;\bm{\lambda})=(1-q^{N-x})(1-p^{-1}q^{-x-1}),\quad
  D(x\,;\bm{\lambda})=p^{-1}q^{N-x}(q^{-x}-1),\\
  &\mathcal{E}(n\,;\bm{\lambda})
  =1-q^n,\quad
  \eta(x\,;\bm{\lambda})=1-q^x,\\
  &P_n(\eta(x\,;\bm{\lambda})\,;\bm{\lambda})
  ={}_2\phi_2\Bigl(
  \genfrac{}{}{0pt}{}{q^{-n},\,q^{-x}}
  {pq,\,q^{-N}}\Bigm|q\,;pq^{x+n+1-N}\Bigr),\\
  &\phi_0(x\,;\bm{\lambda})^2
  =\frac{(q\,;q)_N}{(q\,;q)_x\,(q\,;q)_{N-x}}\,
  \frac{(p^{-1}q^{-N}\,;q)_N}{(p^{-1}q^{-N}\,;q)_{N-x}}\,
  p^xq^{x(x+1-N)},\\
  &d_n(\bm{\lambda})^2
  =\frac{(q\,;q)_N}{(q\,;q)_n\,(q\,;q)_{N-n}}\,
  \frac{(p^{-1}q^{-N}\,;q)_N}{(p^{-1}q^{-N}\,;q)_{N-n}}\,
  p^nq^{n(n+1-N)}
   \times(pq)^{-N},\\
  &A_n(\bm{\lambda})
  =-(1-p^{-1}q^{-n-1})(1-q^{N-n}),\quad
  C_n(\bm{\lambda})
  =-p^{-1}q^{N-n}(q^{-n}-1),\\
  &R_1(z\,;\bm{\lambda})=(q^{-\frac12}-q^{\frac12})^2z',\quad
  z'\eqdef z-1,\\
  &R_0(z\,;\bm{\lambda})=(q^{-\frac12}-q^{\frac12})^2
  z^{\prime\,2},\\
  &R_{-1}(z\,;\bm{\lambda})=(q^{-\frac12}-q^{\frac12})^2
  \Bigl(-z^{\prime\,2}
  -\bigl(p^{-1}q^{-1}+(1+p^{-1})q^N\bigr)z'-p^{-1}(1+q^{-1})q^N\Bigr),\\
  &B(x\,;\bm{\lambda}+\bm{\delta})=B(x+1\,;\bm{\lambda}),\quad
  D(x\,;\bm{\lambda}+\bm{\delta})=q^{-2}D(x\,;\bm{\lambda}),\\
  &\varphi(x\,;\bm{\lambda})=q^x,\quad
  f_n(\bm{\lambda})=\mathcal{E}(n\,;\bm{\lambda}),\quad
  b_n(\bm{\lambda})=1.
\end{align}
The functions $B(x)$ and $D(x)$ are related to those of the quantum
$q$-Krawtchouk polynomial with the change of variable
$x\leftrightarrow N-x$, $B\leftrightarrow D$:
\begin{equation}
  B(x)=D^{\text{quantum $q$-Krawtchouk}}(N-x),\quad
  D(x)=B^{\text{quantum $q$-Krawtchouk}}(N-x).
\end{equation}

\section{Summary and Comments}
\setcounter{equation}{0}

A unified theory of orthogonal polynomials of a discrete variable is
presented through the eigenvalue problem of hermitian matrices of
finite or infinite dimensions.
The hermitian matrices ($\mathcal{H}$) are real and tridiagonal
(Jacobi) matrices \eqref{genham} that can be factorised
$\mathcal{H}=\mathcal{A}^\dagger\mathcal{A}$ and thus are positive
semi-definite.
The orthogonality measure $\phi_0(x)^2$ is obtained as a solution
of the factored equation $\mathcal{A}\phi_0(x)=0$, as in the
ordinary quantum mechanics.
Then the eigenvalue problem can be solved in two different ways;
the first as a difference equation and the second through the three
term recurrence relations.
The former gives the eigenvector polynomial $P_n(\eta(x))$ and
the latter provides the dual polynomial $Q_x(\mathcal{E}(n))$
satisfying the relation $P_n(\eta(x))=Q_x(\mathcal{E}(n))$ through
the universal normalisation condition $P_n(0)=1=Q_x(0)$.
Here $\mathcal{E}(n)$ is the eigenvalue of $\mathcal{H}$ and
$\eta(x)$ is the sinusoidal coordinate satisfying the closure relation
$[\mathcal{H},[\mathcal{H},\eta(x)]]=\eta(x)R_0(\mathcal{H})
+[\mathcal{H},\eta(x)]R_1(\mathcal{H})+R_{-1}(\mathcal{H})$.
Thanks to the closure relation and the shape invariance,
the entire spectrum, the coefficients $A_n$ and $C_n$ of the
three term recurrence relation of $P_n(\eta)$, its normalisation
constants $d_n^2$ are determined algebraically.
As byproducts, some as yet unexplored (not reported in the review
of Koekoek-Swarttouw \cite{koeswart}) dual polynomials with
explicit forms of the difference equation, three term recurrence,
the normalisation measure/constants are presented in
\S\ref{alternatives}.
On top of them we mentioned the dual quantum $q$-Krawtchouk
\eqref{dualqqkrawt1}--\eqref{dualqqkrawt3},
the dual $q$-Meixner \eqref{dualqmeixner1}--\eqref{dualqmeixner3}
and the dual $q$-Charlier \eqref{dualqcharlier1}--\eqref{dualqcharlier3}.

After completing the main part of the paper, we became aware of
the work of Terwilliger \cite{terw} on the correspondence between
the Leonard pair \cite{leonard} and a class of orthogonal polynomials.
Although the general setting of the problem and the methods are
markedly different, Terwilliger's results have some overlap with
the finite dimensional case (I) of the present paper.

\section*{Acknowledgements}

This work is supported in part by Grants-in-Aid for Scientific
Research from the Ministry of Education, Culture, Sports, Science
and Technology, No.18340061 and No.19540179.
This work was also supported in part by the Italian MIUR
(Internazionalizzazione Program) within the joint SISSA--YITP
research project on ``Fundamental Interactions and the Early Universe."

\section*{Appendix A: Possible Forms of the Hamiltonians}
\label{appendA}
\setcounter{equation}{0}
\renewcommand{\theequation}{A.\arabic{equation}}

Here we continue from \S\ref{detetaBD} the algebraic analysis
of the closure relation to determine the possible function forms
of $B(x)$ and $D(x)$. The remaining equation to be solved is
\eqref{crcond3}, which is the diagonal part of the closure relation
\eqref{closurerel1t}.

With \eqref{cransr} and \eqref{crcond1pp}, the equation to be
analysed \eqref{crcond3} reads
\begin{align}
  &\bigl(\eta(x)-\eta(x+2)\bigr)B(x)D(x+1)
  +\bigl(\eta(x)-\eta(x-2)\bigr)B(x-1)D(x)\n
  &\qquad\qquad
  =\bigl(\eta(x+1)-2\eta(x)+\eta(x-1)\bigr)a_x^2
  +\bigl(2r_1^{(0)}\eta(x)+r_{-1}^{(1)}\bigr)a_x
  +r_0^{(0)}\eta(x)+r_{-1}^{(0)}.
  \label{crcond3pp}
\end{align}
Before proceeding, let us stress again that the input of the
eigenvalue problem \eqref{eigenpr} is $B(x)$ and $D(x)$.
The quantities $r_i^{(j)}$ ($i=1,0,-1$; $j=0,1,2$) are determined
by the consistency of the closure relation \eqref{closurerel1}.
In Appendix A, however, we turn the logic and regard $r_i^{(j)}$ as
input data and determine the possible forms of $B(x)$ and $D(x)$
which satisfy \eqref{crcond1}--\eqref{crcond3}.
In this approach the positivity of $B(x)$ and $D(x)$ must be verified
at the last step, since it is not guaranteed automatically.
By multiplying $\eta(x-1)-\eta(x+1)$ to \eqref{crcond3pp}, we obtain
\begin{align}
  &\ \
  \bigl(\eta(x+2)-\eta(x)\bigr)\bigl(\eta(x+1)-\eta(x-1)\bigr)B(x)D(x+1)
  -(x\to x-1)\n
  &=\Bigl(\bigl(\eta(x-1)-\eta(x)\bigr)^2
  -\bigl(\eta(x+1)-\eta(x)\bigr)^2\Bigr)a_x^2\n
  &\quad
  +\bigl(\eta(x-1)-\eta(x+1)\bigr)
  \bigl(2r_1^{(0)}\eta(x)+r_{-1}^{(1)}\bigr)a_x
  +\bigl(\eta(x-1)-\eta(x+1)\bigr)\bigl(r_0^{(0)}\eta(x)+r_{-1}^{(0)}\bigr)\n
  &=
  \frac{\bigl(r_1^{(0)}\eta(x+1)\eta(x)+r_{-1}^{(1)}\eta(x+1)+a'_0\bigr)
        \bigl(r_1^{(0)}\eta(x+1)\eta(x)+r_{-1}^{(1)}\eta(x)+a'_0\bigr)}
       {\bigl(\eta(x+1)-\eta(x)\bigr)^2}\n
  &\qquad
  -r_0^{(0)}\eta(x+1)\eta(x)
  -r_{-1}^{(0)}\bigl(\eta(x+1)+\eta(x)\bigr)
  -(x\to x-1),
\end{align}
where $a'_0=\eta(1)\eta(-1)a_0$ and $a_0=B(0)$.
For the second equality, the main result of \S\ref{detetaBD},
\eqref{ansa} is used. The above equation is simplified as
\begin{align}
  &\bigl(\eta(x+2)-\eta(x)\bigr)\bigl(\eta(x+1)-\eta(x-1)\bigr)B(x)D(x+1)\n
  &=
  \frac{\bigl(r_1^{(0)}\eta(x+1)\eta(x)+r_{-1}^{(1)}\eta(x+1)+a'_0\bigr)
        \bigl(r_1^{(0)}\eta(x+1)\eta(x)+r_{-1}^{(1)}\eta(x)+a'_0\bigr)}
       {\bigl(\eta(x+1)-\eta(x)\bigr)^2}\n
  &\quad
  -\eta(1)a_0\bigl(r_{-1}^{(1)}+\eta(1)a_0\bigr)
  -r_0^{(0)}\eta(x+1)\eta(x)
  -r_{-1}^{(0)}\bigl(\eta(x+1)+\eta(x)-\eta(-1)\bigr).
  \label{ansBD}
\end{align}

For each case of $\eta(x)$, \eqref{etaform1}--\eqref{etaform5},
we will solve $B(x)$ from the above equation \eqref{ansBD}, with
the result of $a_x=B(x)+D(x)$, \eqref{ansa}.
Let us parametrise $B(x)$ as
\begin{equation}
  B(x)=\frac{\widetilde{B}(x)}
  {\bigl(\eta(x+1)-\eta(x)\bigr)\bigl(\eta(x+1)-\eta(x-1)\bigr)}.
  \label{B=Bt/etaeta}
\end{equation}
Note that the denominator is the dual counterpart of
$\alpha_+(\mathcal{E}(n))\bigl(\alpha_+(\mathcal{E}(n))-
\alpha_-(\mathcal{E}(n))\bigr)$
in the denominator of $A_n$ in \eqref{Anformula2}.

We eliminate $D(x+1)$ from \eqref{ansBD} (multiplied by
$\bigl(\eta(x+2)-\eta(x+1)\bigr)\bigl(\eta(x+1)-\eta(x)\bigr)^2$)
by using \eqref{ansa}:
\begin{align}
  &-\widetilde{B}(x)\left(\!
  \bigl(r_1^{(0)}\eta(x+1)^2+r_{-1}^{(1)}\eta(x+1)+a'_0\bigr)
  \bigl(\eta(x+2)-\eta(x)\bigr)
  \!+\widetilde{B}(x+1)\bigl(\eta(x+1)-\eta(x)\bigr)\!\right)\n
  &=\bigl(\eta(x+2)-\eta(x+1)\bigr)\n
  &\quad\times\biggl(
  \bigl(r_1^{(0)}\eta(x+1)\eta(x)+r_{-1}^{(1)}\eta(x+1)+a'_0\bigr)
  \bigl(r_1^{(0)}\eta(x+1)\eta(x)+r_{-1}^{(1)}\eta(x)+a'_0\bigr)\n
  &\qquad\quad
  -\bigl(\eta(x+1)-\eta(x)\bigr)^2\Bigl(
  \eta(1)a_0\bigl(r_{-1}^{(1)}+\eta(1)a_0\bigr)
  +r_0^{(0)}\eta(x+1)\eta(x)\n
  &\qquad\qquad\qquad\qquad\qquad\qquad\qquad\qquad
  +r_{-1}^{(0)}\bigl(\eta(x+1)+\eta(x)-\eta(-1)\bigr)\Bigr)\biggr).
  \label{ansBD2}
\end{align}
This equation determines $\widetilde{B}(x)$ uniquely from the
boundary (initial) value
$\widetilde{B}(0)=\eta(1)\bigl(\eta(1)$ $-\eta(-1)\bigr)B(0)$.

We have carried out this last step by assuming that $\widetilde{B}(x)$
is a polynomial in $x$, $\widetilde{B}(x)=\sum_{j=0}^M\beta_jx^j$
for the cases (\romannumeral1) and (\romannumeral2) and that
$\widetilde{B}(x)$ is a Laurent polynomial in $q^x$,
$\widetilde{B}(x)=\sum_{j=M_1}^{M_2}\beta_jq^{jx}$
for the cases (\romannumeral3) -- (\romannumeral5).
The results are as follows.

\noindent
\underline{case (\romannumeral1),\quad $\eta(x)=x$}:\\
A quadratic polynomial ($M=2$) solution  $\widetilde{B}(x)$ is
obtained and $B(x)$ and $D(x)$ are:
\begin{equation}
  B(x)=(\text{quadratic in $x$}),\quad D(x)=x(\text{linear in $x$}),
\end{equation}
which contains the Hahn polynomial \S\ref{[KS1.5]} as the most
generic case.

\noindent
\underline{case (\romannumeral2),\quad $\eta(x)=\epsilon'x(x+d)$}:\\
A quartic polynomial ($M=4$) solution $\widetilde{B}(x)$ is obtained
and $B(x)$ and $D(x)$ are:
\begin{equation}
  B(x)=\frac{(x+d)(\text{cubic in $x$})}{(2x+d)(2x+1+d)},\quad
  D(x)=\frac{x(\text{cubic in $x$})}{(2x+d)(2x-1+d)},
\end{equation}
which contains the Racah polynomial \S\ref{[KS1.2]} as the most
generic case.

\noindent
\underline{case (\romannumeral3),\quad $\eta(x)=1-q^x$}:\\
A quadratic polynomial ($M_1=0$, $M_2=2$) solution $\widetilde{B}(x)$
in $q^x$ is obtained and $B(x)$ and $D(x)$ are:
\begin{equation}
  B(x)=(\text{quadratic in $q^{-x}$}),\quad
  D(x)=(q^{-x}-1)(\text{linear in $q^{-x}$}).
  \label{BDq3}
\end{equation}
The most generic one in this category reported in the
Koekoek-Swarttouw's review \cite{koeswart} is the little $q$-Jacobi
polynomial \S\ref{[KS3.12]}. The dual quantum $q$-Krawtchouk polynomial
\eqref{dualqqkrawt1}--\eqref{dualqqkrawt3} also belongs to this case.
There is a more general polynomial with one more parameter as
tentatively called alternative $q$-Hahn polynomial \S\ref{aqHahn}.
For the infinite dimensional case (II), the functions $B(x)$ and
$D(x)$ grow exponentially. Therefore the corresponding spectrum is
unbounded, {\em i.e.\/} type (\romannumeral4) \eqref{matheform4} or
(\romannumeral5) \eqref{matheform5}. In other words, there is no
self-dual $\eta(x)=1-q^x$, $\mathcal{E}(n)=1-q^n$ polynomial in
the infinite dimensional case (II).

\noindent
\underline{case (\romannumeral4),\quad $\eta(x)=q^{-x}-1$}:\\
A quadratic polynomial ($M_1=-2$, $M_2=0$) solution $\widetilde{B}(x)$
in $q^{-x}$ is obtained and $B(x)$ and $D(x)$ are:
\begin{equation}
  B(x)=(\text{quadratic in $q^x$}),\quad
  D(x)=(1-q^{x})(\text{linear in $q^x$}),
\end{equation}
which contains the $q$-Hahn polynomial \S\ref{[KS3.6]} as the most
generic one.
In contrast to the previous case the functions $B(x)$ and $D(x)$
are all bounded in the the infinite dimensional case (II).
Therefore the corresponding spectrum is also bounded, {\em i.e.\/}
type (\romannumeral3) \eqref{matheform3} only.
In other words, there is no self-dual $\eta(x)=q^{-x}-1$,
$\mathcal{E}(n)=q^{-n}-1$ polynomial in the infinite dimensional
case (II).

\noindent
\underline{case (\romannumeral5),\quad
$\eta(x)=\epsilon'(q^{-x}-1)(1-dq^x)$}:\\
A quadratic polynomial ($M_1=-2$, $M_2=2$) solution $\widetilde{B}(x)$
in $q^{x}$ and $q^{-x}$ is obtained and $B(x)$ and $D(x)$ are:
\begin{equation}
  B(x)=\frac{(1-dq^{x})(\text{cubic in $q^x$})}
       {(1-dq^{2x})(1-dq^{2x+1})},\quad
  D(x)=\frac{(1-q^x)(\text{cubic in $q^x$})}{(1-dq^{2x})(1-dq^{2x-1})},
\end{equation}
which contains the most generic case of the $q$-Racah polynomial
\S\ref{[KS3.2]}. The functions $B(x)$ and $D(x)$ are all bounded
in the the infinite dimensional case (II).
The corresponding spectrum is also bounded, {\em i.e.\/}
type (\romannumeral3) \eqref{matheform3} only.
The case (\romannumeral5) does not have self-dual polynomial, either.
Thus the most generic functional forms of $B(x)$ and $D(x)$ are
determined for each case of the five different sinusoidal coordinates
(\romannumeral1)--(\romannumeral5).

Determination of polynomials satisfying certain forms of difference
equations for given $\eta(x)$ (quadratic or $q$-quadratic) has a
long history \cite{askeywil,vinzhed,leonard,terw}.
Bochner's theorem \cite{bochner} on Sturm-Liouville polynomials is
a precursor.
The present characterisation in terms of the closure relation
\eqref{closurerel1t} is consistent with the existing ones.

\bigskip

It is interesting to note that in all the known cases the following
relation is satisfied:
\begin{equation}
  \frac{r_0^{(0)}}{B(0)^2}
  +\frac{r_1^{(0)}}{B(0)}\frac{r_{-1}^{(0)}}{\eta(1)B(0)^2}
  -\Bigl(\frac{r_{-1}^{(0)}}{\eta(1)B(0)^2}\Bigr)^2=0.
  \label{ri0cond}
\end{equation}
When this condition is met, the solution of \eqref{ansBD2} for all
the cases (\romannumeral1)--(\romannumeral5) has a simple expression
in terms of $\eta(x)$:
\begin{align}
  \widetilde{B}(x)&=-r_1^{(0)}\eta(x)\eta(x+1)
  +\frac{r_{-1}^{(0)}}{\eta(1)B(0)}\eta(x)\bigl(\eta(x+1)-\eta(x)\bigr) \n
  &\quad
  -r_{-1}^{(1)}\eta(x)+\eta(1)B(0)\bigl(\eta(x+1)-\eta(x)-\eta(-1)\bigr).
  \label{Btilderes}
\end{align}
For the other function $D(x)$ we parametrise
\begin{equation}
  D(x)=\frac{\widetilde{D}(x)}
  {\bigl(\eta(x-1)-\eta(x)\bigr)\bigl(\eta(x-1)-\eta(x+1)\bigr)}.
  \label{D=Dt/etaeta}
\end{equation}
The corresponding expression for $\widetilde{D}(x)$ is
\begin{align}
  \widetilde{D}(x)&=-r_1^{(0)}\eta(x)\eta(x-1)
  +\frac{r_{-1}^{(0)}}{\eta(1)B(0)}\eta(x)\bigl(\eta(x-1)-\eta(x)\bigr) \n
  &\quad
  -r_{-1}^{(1)}\eta(x)+\eta(1)B(0)\bigl(\eta(x-1)-\eta(x)-\eta(-1)\bigr).
  \label{Dtilderes}
\end{align}

The above relation \eqref{ri0cond} is a consequence of \eqref{ABrel},
which can be regarded as a consistency condition.
{}From \eqref{alphapmexp}, \eqref{alphapmE} and \eqref{A0exp} we have
$\mathcal{E}(1)=\alpha_+\bigl(\mathcal{E}(0))=\frac12(r_1^{(0)}
+\sqrt{r_1^{(0)\,2}+4r_0^{(0)}}\,\bigr)$ and
$A_0=r_{-1}^{(0)}/r_0^{(0)}$. Substituting these into \eqref{ABrel}
(multiplied by $r_0^{(0)}/\mathcal{E}(1)$),
we obtain
\[
  r_{-1}^{(0)}+\frac12\Bigl(\sqrt{r_1^{(0)\,2}+4r_0^{(0)}}-r_1^{(0)}\Bigr)
  \eta(1)B(0)=0,
\]
namely
\[
  \frac12\sqrt{r_1^{(0)\,2}+4r_0^{(0)}}
  =\frac12r_1^{(0)}-\frac{r_{-1}^{(0)}}{\eta(1)B(0)}.
\]
Squaring this gives \eqref{ri0cond}.

By using \eqref{ABrel}, the expressions for $B(x)$ and $D(x)$,
\eqref{B=Bt/etaeta} with \eqref{Btilderes} and
\eqref{D=Dt/etaeta} with \eqref{Dtilderes}, are cast
in a similar form as $A_n$ \eqref{Anformula2} and $C_n$
\eqref{Cnformula2}:
\begin{align}
  B(x)&=-\frac{R^{\text{dual}}_{-1}(\eta(x))
  +\mathcal{E}(1)\bigl(\eta(x)+A_0\bigr)\alpha^{\text{dual}}_+(\eta(x))}
  {\alpha^{\text{dual}}_+(\eta(x))
  \bigl(\alpha^{\text{dual}}_+(\eta(x))
  -\alpha^{\text{dual}}_-(\eta(x))\bigr)}\,,
  \label{Bxformula2}\\
  D(x)&=-\frac{R^{\text{dual}}_{-1}(\eta(x))
  +\mathcal{E}(1)\bigl(\eta(x)+A_0\bigr)\alpha^{\text{dual}}_-(\eta(x))}
  {\alpha^{\text{dual}}_-(\eta(x))
  \bigl(\alpha^{\text{dual}}_-(\eta(x))
  -\alpha^{\text{dual}}_+(\eta(x))\bigr)}\,.
  \label{Dxformula2}
\end{align}
Here $R^{\text{dual}}_i(z)$ and $\alpha^{\text{dual}}_{\pm}(z)$ are
functions appearing in the dual closure relation
\S\ref{dualclosuresection}.
The functions $R^{\text{dual}}_i(z)$ are given in 
\eqref{dualclR1}--\eqref{dualclR-1} and $\alpha^{\text{dual}}_{\pm}(z)$
are defined by (c.f. \eqref{alphapmexp}, \eqref{alphapmE})
\begin{equation}
  \alpha^{\text{dual}}_{\pm}(z)\eqdef\frac12\Bigl(R^{\text{dual}}_1(z)
  \pm\sqrt{R^{\text{dual}}_1(z)^2+4R^{\text{dual}}_0(z)}\,\Bigr),
\end{equation}
which satisfy
\begin{equation}
  \alpha^{\text{dual}}_{\pm}(\eta(x))=\eta(x\pm 1)-\eta(x).
\end{equation}
Similarly $A_n$ and $C_n$ can be cast into the forms like $B(x)$ in
\eqref{B=Bt/etaeta} with \eqref{Btilderes} and $D(x)$ in
\eqref{D=Dt/etaeta} with \eqref{Dtilderes}.

\section*{Appendix B: Some definitions related to the hypergeometric
and $q$-hypergeometric functions}
\label{appendB}
\setcounter{equation}{0}
\renewcommand{\theequation}{B.\arabic{equation}}

For self-containedness we collect several definitions related to
the ($q$-)hypergeometric functions \cite{koeswart}.

\noindent
$\circ$ Pochhammer symbol $(a)_n$ :
\begin{equation}
  (a)_n\eqdef\prod_{k=1}^n(a+k-1)=a(a+1)\cdots(a+n-1)
  =\frac{\Gamma(a+n)}{\Gamma(a)}.
  \label{defPoch}
\end{equation}
$\circ$ $q$-Pochhammer symbol $(a\,;q)_n$ :
\begin{equation}
  (a\,;q)_n\eqdef\prod_{k=1}^n(1-aq^{k-1})=(1-a)(1-aq)\cdots(1-aq^{n-1}).
  \label{defqPoch}
\end{equation}
$\circ$ hypergeometric series ${}_rF_s$ :
\begin{equation}
  {}_rF_s\Bigl(\genfrac{}{}{0pt}{}{a_1,\,\cdots,a_r}{b_1,\,\cdots,b_s}
  \Bigm|z\Bigr)
  \eqdef\sum_{n=0}^{\infty}
  \frac{(a_1,\,\cdots,a_r)_n}{(b_1,\,\cdots,b_s)_n}\frac{z^n}{n!}\,,
  \label{defhypergeom}
\end{equation}
where $(a_1,\,\cdots,a_r)_n\eqdef\prod_{j=1}^r(a_j)_n
=(a_1)_n\cdots(a_r)_n$.\\
$\circ$ $q$-hypergeometric series (the basic hypergeometric series)
${}_r\phi_s$ :
\begin{equation}
  {}_r\phi_s\Bigl(
  \genfrac{}{}{0pt}{}{a_1,\,\cdots,a_r}{b_1,\,\cdots,b_s}
  \Bigm|q\,;z\Bigr)
  \eqdef\sum_{n=0}^{\infty}
  \frac{(a_1,\,\cdots,a_r\,;q)_n}{(b_1,\,\cdots,b_s\,;q)_n}
  (-1)^{(1+s-r)n}q^{(1+s-r)n(n-1)/2}\frac{z^n}{(q\,;q)_n}\,,
  \label{defqhypergeom}
\end{equation}
where $(a_1,\,\cdots,a_r\,;q)_n\eqdef\prod_{j=1}^r(a_j\,;q)_n
=(a_1\,;q)_n\cdots(a_r\,;q)_n$.


\end{document}